\documentclass[11pt, reqno]{amsart}
\usepackage{amssymb}
\usepackage{eucal}
\usepackage{amsmath}
\usepackage{amscd}
\usepackage[dvips]{color}
\usepackage{multicol}
\usepackage[all]{xy}           
\usepackage{graphicx}
\usepackage{color}
\usepackage{colordvi}
\usepackage{xspace}
\usepackage{tikz}
\usepackage{cite}
\usepackage{graphicx}
\usepackage{subfigure}
\graphicspath{{Figures/}}
\usepackage{ifpdf}
\ifpdf
 \usepackage[colorlinks,final,hyperindex]{hyperref}
\else
 \usepackage[colorlinks,final,hyperindex,hypertex]{hyperref}
\fi
\hypersetup{colorlinks=true,linkcolor=black}

\begin{document}
\newcommand {\emptycomment}[1]{} 

\baselineskip=14pt
\newcommand{\nc}{\newcommand}
\newcommand{\delete}[1]{}
\nc{\mfootnote}[1]{\footnote{#1}} 
\nc{\todo}[1]{\tred{To do:} #1}

\def\a{\alpha}
\def\b{\beta}
\def\bd{\boxdot}
\def\bbf{\bar{f}}
\def\bF{\bar{F}}
\def\bbF{\bar{\bar{F}}}
\def\bbbf{\bar{\bar{f}}}
\def\bg{\bar{g}}
\def\bG{\bar{G}}
\def\bbG{\bar{\bar{G}}}
\def\bbg{\bar{\bar{g}}}
\def\bT{\bar{T}}
\def\bt{\bar{t}}
\def\bbT{\bar{\bar{T}}}
\def\bbt{\bar{\bar{t}}}
\def\bR{\bar{R}}
\def\br{\bar{r}}
\def\bbR{\bar{\bar{R}}}
\def\bbr{\bar{\bar{r}}}
\def\bu{\bar{u}}
\def\bU{\bar{U}}
\def\bbU{\bar{\bar{U}}}
\def\bbu{\bar{\bar{u}}}
\def\bw{\bar{w}}
\def\bW{\bar{W}}
\def\bbW{\bar{\bar{W}}}
\def\bbw{\bar{\bar{w}}}
\def\btl{\blacktriangleright}
\def\btr{\blacktriangleleft}
\def\c{\cdot}
\def\ci{\circ}
\def\d{\delta}
\def\dd{\diamondsuit}
\def\D{\Delta}
\def\G{\Gamma}
\def\g{\gamma}
\def\l{\lambda}
\def\lr{\longrightarrow}
\def\o{\otimes}
\def\om{\omega}
\def\p{\psi}
\def\r{\rho}
\def\ra{\rightarrow}
\def\rh{\rightharpoonup}
\def\lh{\leftharpoonup}
\def\s{\sigma}
\def\st{\star}
\def\ti{\times}
\def\tl{\triangleright}
\def\tr{\triangleleft}
\def\v{\varepsilon}
\def\vp{\varphi}

\newtheorem{thm}{Theorem}[section]
\newtheorem{lem}[thm]{Lemma}
\newtheorem{cor}[thm]{Corollary}
\newtheorem{pro}[thm]{Proposition}
\theoremstyle{definition}
\newtheorem{defi}[thm]{Definition}
\newtheorem{ex}[thm]{Example}
\newtheorem{rmk}[thm]{Remark}
\newtheorem{pdef}[thm]{Proposition-Definition}
\newtheorem{condition}[thm]{Condition}
\newtheorem{question}[thm]{Question}
\renewcommand{\labelenumi}{{\rm(\alph{enumi})}}
\renewcommand{\theenumi}{\alph{enumi}}

\nc{\ts}[1]{\textcolor{purple}{Tianshui:#1}}

\font\cyr=wncyr10

 \title{\bf An associative analogy of Lie $H$-pseudobialgebra}

 \author{Linlin Liu\textsuperscript{*}}
 \address{Department of Science, Henan Institute of Technology, Xinxiang 453003, China}
         \email{liulinlin2016@163.com}

  \author{Zhitao Guo}
 \address{Department of Science, Henan Institute of Technology, Xinxiang 453003, China}
         \email{guotao60698@163.com}

  \thanks{\textsuperscript{*}Corresponding author}

\date{\today}

\begin{abstract}
The purpose of this paper is to study infinitesimal $H$-pseudobialgebra, which is an associative analogy of Lie $H$-pseudobialgebra. We first define the infinitesimal $H$-pseudobialgebra and investigate some properties of this new algebraic structure. Then we consider the coboundary infinitesimal $H$-pseudobialgebra, which is the subclass of infinitesimal $H$-pseudobialgebra and we obtain the associative Yang-Baxter equation over an associative $H$-pseudoalgebra. Finally, we found the connection between the (coboundary) infinitesimal $H$-pseudobialgebra and the (coboundary) Lie $H$-pseudobialgebra. Meanwhile, the relationship between the associative Yang-Baxter equation and the classical Yang-Baxter equation (over an $H$-pseudoalgebra) is established.
\end{abstract}



 \maketitle

\tableofcontents

\allowdisplaybreaks

\section{Introduction}

The notion of conformal algebra (\cite{K1}) was introduced by Kac as an axiomatic description of the operator product expansion (OPE) of chiral fields in conformal field theory, and it came to be useful for investigation of vertex algebras.  Recall that a Lie conformal algebra $L$ is defined as a $\mathbb{C}[\partial]$-module ($\partial$ is an indeterminate), endowed with a $\mathbb{C}$-linear map
$$
L\o L\longrightarrow \mathbb{C}[\lambda]\o L,\quad a\o b\mapsto [a_{\lambda}b]
$$
satisfying axioms similar to those of Lie algebra (see \cite{DK, K1}). Later, Bakalov, D'Andrea and Kac replaced the above polynomial algebra $\mathbb{C}[\partial]$ with any cocommutative Hopf algebra $H$ in \cite{BDK}, and found that this high-dimensional conformal algebra is actually the algebra in a pseudotensor category, which is called (Lie) $H$-pseudoalgebra. So far, the classification problems, cohomology theory and representation theory of $H$-pseudoalgebras have been considered in \cite{BDK, BDK2, BKV}.
As a natural generation of conformal algebras, (Lie) $H$-pseudoalgebras are closely related to the differential Lie algebras of the Ritt and Hamiltonian formalism in the theory of nonlinear evolution equations (see \cite{D, GD1, GD2}). However, their roles in many fields of mathematical physics are not yet completely understood since they are relatively new algebraic structures. The fact that the annihilation algebra of the associative $H$-pseudoalgebra Cend$H$ is nothing else but the Drinfeld Double of the Hopf algebra $H$, which leads us to believe that there should be a deep connection between the theory of $H$-pseudoalgebras and quantum groups. Recently there has been some interest in the theory of $H$-pseudoalgebras (see, for example, \cite{BL, LW, S1, SW, WZ}).

Infinitesimal bialgebras (also called $\v$-bialgebras) were introduced by Join and Rota in order to provide an algebraic framework for the calculus of
divided differences (\cite{JR}). More precisely, an infinitesimal bialgebra $(A, \mu,\D)$ is an associative algebra $(A,\mu)$ and a coassociative coalgebra $(A,\D)$ such that the comultiplication $\D$ is a 1-cocycle in algebra cohomology (i.e., a derivation) with coefficient in $A\o A$. Moreover, the notions of coboundary and quasitriangular infinitesimal bialgebra, infinitesimal Hopf algebra, and the basic theory were established by Aguiar in \cite{A4, A1}. Further research can be found in \cite{A3, WW}.

In addition to the widely research in combinatorics (\cite{A2, ER, HR}), infinitesimal bialgebras are also closely related to Lie bialgebras (\cite{D1, D2, JR}). Indeed, a Lie bialgebra is a Lie algebra and a Lie coalgebra, in which the cobracket is a 1-cocycle in Lie algebra cohomology. Thus the cocycle condition in an infinitesimal bialgebra can be seen as an associative analogy of that in a Lie bialgebra. The necessary and sufficient conditions for infinitesimal bialgebra to be a related Lie bialgebra were given in \cite{A1}. Furthermore, it was shown that the solutions of associative Yang-Baxter equation (\cite{A1, Ho}) can induce the solutions of classical Yang-Baxter equation (\cite{BGN, Bou}).

Lie $H$-pseudobialgebra appeared in \cite{BL} as a generalization of conformal bialgebra (\cite{L}) and Lie bialgebra. Compared with Lie bialgebra, Lie $H$-pseudobialgebra is defined on the $H$-modules rather than the vector spaces. More precisely, Lie $H$-pseudobialgebra is both a Lie $H$-coalgebra and a Lie $H$-pseudoalgebra satisfying the cocycle condition. Similar to the relationship between infinitesimal bialgebra and Lie bialgebra, we want to find an associative analog of Lie $H$-pseudobialgebra, which is our motivation for defining infinitesimal $H$-pseudobialgebra. In addition, there is a natural, but not obvious way to establish the connections between infinitesimal $H$-pseudobialgebra and Lie $H$-pseudobialgebra. The current paper will devote to these questions.

The paper is organized as follows.

In Section 2, we recall some basic notions about associative $H$-pseudoalgebra and present some properties of its representation.

In Section 3, we mainly define the notion of infinitesimal $H$-pseudobialgebra (see Definition 3.4). Specifically, an infinitesimal $H$-pseudobialgebra $(A, *,\D)$ consists of an associative $H$-pseudoalgebra $(A,*)$ and a coassociative $H$-coalgebra $(A,\D)$ such that $\D$ is a 1-cocycle in associative $H$-pseudoalgebra cohomology (see Section 4) with coefficient in $A^{\o 2}$, which is an associative analoy of the cocycle condition in a Lie $H$-pseudobialgebra. Moreover, some examples and basic properties are given.

In Section 4, we study the cohomology theory of associative $H$-pseudoalgebra and the coboundary innitesimal $H$-pseudobialgebra (see Definition 4.2), which is an important subclass of infinitesimal $H$-pseudobialgebra. In particular, we consider the construction of coboundary infinitesimal $H$-pseudobialgebra and the corresponding associative pseudo-Yang-Baxter equation ($pseudo$-AYBE)(see Theorem 4.4).

In Sections 5, we mainly study the relationship between infinitesimal $H$-pseudobialgebras and Lie $H$-pseudobialgebras (see Theorem 5.6). We first find the necessary and sufficient conditions for infinitesimal $H$-pseudobialgebra to be related Lie $H$-pseudobialgebra (see Corollary 5.7), then we give a sufficient condition under which a coboundary infinitesimal $H$-pseudobialgebra gives rise to the corresponding coboundary Lie $H$-pseudobialgebra (see Theorem 5.10).

In the last sections, we show that under some suitable conditions, a solution of the $pseudo$-AYBE is also a solution of the classical Yang-Baxter equation ($pseudo$-CYBE) in related Lie $H$-pseudobialgebra (see Theorem 6.1).

Throughout this paper, $k$ is a fixed algebraically closed field of characteristic zero, $H$ is a cocommutative Hopf algebra over $k$, and $X=H^{*}=Hom_{k}(H, k)$ denotes the dual of $H$.  As usual, we adopt Sweedler's notations in \cite{S}. For a coalgebra $C$, we write its comultiplication as $\D(c)=c_1\o c_2$, $\forall~c\in C$. For any vector space $V$, we will define $\sigma(f\o g)=g\o f$, $(12)(f\o g\o h)=g\o f\o h$ and $(13)(f\o g\o h)=h\o g\o f$, for all $f, g, h\in V$.

\section{Associative $H$-pseudoalgebras and their representations}\def\theequation{2.\arabic{equation}}
\setcounter{equation} {0}

We first recall some definitions and notations of associative $H$-pseudoalgebras (see \cite{BDK, BL} for more details).
\\

{\bf Definition 2.1} An $H$-pseudoalgebra is a left $H$-module $A$ together with a map (called the pseudoproduct):
$$
*: \quad A\o A\longrightarrow (H\o H)\o_H A,~~~a\o b\mapsto a*b
$$
satisfying
{\bf $H$-bilinearity:} \quad for all $a, b\in A$ and $f, g\in H$, one has
$$
fa*gb=(£¨f\o g£©\o_H 1)(a*b).
$$
If $a*b=\sum_{i}(f_i\o g_i)\o_H e_i$, then we have $fa*gb=\sum_{i}(ff_i\o gg_i)\o_H e_i$.

Note that the $H$-pseudoalgebra $(A, *)$ is associative if it satisfies $(a* b)* c=a*(b*c)$ in $H^{\o 3}\o_H A$, for any
 $a, b, c\in A$, and $(A, *)$ is commutative if $a*b=(\sigma\o_H id)(b*a)$ holds.

Moreover, $H$-pseudoalgebra $(A, *)$ is called finite, if it is finitely generated as $H$-module.

In particular, for the one-dimensional Hopf algebra $H=k$, an associative $H$-pseudoalgebra is just an ordinary algebra over the field $k$.
\\

For an arbitrary Hopf algebra $H$, we recall that the map $\mathcal{F}: H\o H\lr H\o H$ defined by the formula
$$
\mathcal{F}(f\o g)=(f\o 1)(S\o id)\D(g)=fS(g_1)\o g_2
$$
is called the Fourier transform. Observe that $\mathcal{F}$ is a vector space isomorphism with an inverse given by
$$
\mathcal{F}^{-1}(f\o g)=(f\o 1)\D(g)=fg_1\o g_2.
$$
In order to reformulate the definition of associative $H$-pseudoalgebra, the authors in \cite{BDK} introduced another product $\{a, b\}\in H\o A$ as the Fourier transform of $a*b$:
\begin{eqnarray*}
&&\{a, b\}=\sum_{i}\mathcal{F}(f_i\o g_i)(1\o e_i)=\sum_{i}f_iS(g_{i1})\o g_{i2}e_i, \hbox{~if~}a*b=\sum_{i} (f_i\o g_i)\o_H e_i.
\end{eqnarray*}
In other words,
$$
\{a, b\}=\sum_{i} h_i\o c_i,~~\hbox{if~}a*b=\sum_{i}(h_i\o 1)\o_H c_i.
$$
For $x\in X:=H^{*}$, the $x$-product in $A$ is as follows:
$$
a\circ_{x}b=(<S(x), \cdot>\o id)[a, b]=\sum_{i}<S(x), f_iS(g_{i1})>g_{i2}e_i=\sum_{i}<S(x), h_i>c_i,
$$
if $a*b=\sum_{i}(f_i\o g_i)\o_H e_i=\sum_{i}(h_i\o 1)\o_H c_i$.

Using properties of the Fourier transform, an equivalent definition of an associative $H$-pseudoalgebra is as follows.
\\

{\bf Definition 2.2} An associative $H$-conformal algebra is a left $H$-module $A$ endowed with a product $\{\cdot, \cdot\}: A\o A\lr H\o A$, satisfying the following properties (for all $a, b, c\in A$ and $h\in H$):

{\bf $H$-sesquilinearity:}
$$
\{ha, b\}=(h\o 1)\{a,b\}, \quad \{a, hb\}=(1\o h_2)\{a, b\}(S(h_1)\o 1).
$$

{\bf Associativity:}
$$
\{a, \{b, c\}\}=(\mathcal{F}^{-1}\o id)\{\{a, b\}, c\}
$$
in $H\o H\o A$, where
$
\{a, \{b, c\}\}=(\sigma\o id)(id\o \{a, \cdot\})\{b, c\}, \{\{a, b\}, c\}=(id\o \{\cdot, c\})\{a, b\}.
$

Definition 2.2 can also be reformulated in terms of the $x$-product. Formally, we will use the same notation for $X$ as for $H$.
\\

{\bf Definition 2.3} An associative $H$-conformal algebra is a left $H$-module $A$ equipped with $x$-product $\circ_{x}: A\o A\longrightarrow A, a\o b\mapsto a\circ_x b$ for all $a, b\in A$ and $x\in X$, satisfying the following properties:

{\bf Locality:} \quad for any basis $\{x_{i}\}$ of $X$, $a\circ_{x_{i}}b\neq 0$ for only a finite number of $i$.

{\bf $H$-sesquilinearity:} $(ha)\circ_xb=a\circ_{xh}b, \quad a\circ_{x}(hb)=h_2(a\circ_{S(h_1)x}b)$, for all $a, b\in A$ and $h\in H$.

{\bf Associativity:} $a\circ_x(b\circ_yc)=(a\circ_{x_2}b)\circ_{yx_1}c$, for all $a, b, c\in A$ and $x, y\in X$.
\\

Let $A, B$ be two associative $H$-pseudoalgebras. A left $A$-module is a left $H$-module $M$ together with an operation $\r_{l}\in Hom_{H\o H}(A\o M, (H\o H)\o_H M)$, we denote $\r_l(a\o m)=a*m$, which satisfies $a*(b* m)=(a* b)* m$, for all $a, b\in A$ and $m\in M$. Similarly, we can define a right $B$-module. If $M$ is a left $A$-module and right $B$-module such that $(a*m)*b=a*(m*b)$ for all $a\in A, b\in B$ and $m\in M$, then $M$ is called an $A$-$B$-bimodule. An $A$-$A$-bimodule is simply called an $A$-bimodule.

Next we show some properties that will be used later. Prior to this, we introduce the following notations:

Suppose that $A$ is an associative $H$-pseudoalgebra and $M$ is a left $A$-module. For all $a\in A$ and $m\in M$, we have $a*m=\sum_{i}f_i\o g_i\o_H m_i=\sum_{i}f_iS(g_{i1})\o 1\o_H g_{i2}m_i\in (H\o H)\o_H M$. By Lemma 2.3 in \cite{BDK}, $a*m$ can be written uniquely in the form $\sum_{i}(h_i\o 1)\o_H c_i$, where $\{h_i\}$ is a fixed $k$-basis of $H$. Throughout this paper, we write
$$
a*m=\sum_{(a, m)}h^{a, m}\o 1\o_H c_{a, m}=h^{a, m}\o 1\o_H c_{a, m}
$$
for convenience. Similarly, for a right $A$-module $N$, we set
$$
n*a=\sum_{(n, a)}1\o l^{n, a}\o_H e_{n, a}=1\o l^{n, a}\o_H e_{n, a}, \forall~a\in A, n\in N.
$$

{\bf Proposition 2.4} Let $(A,*)$ be an associative $H$-pseudoalgebra. Suppose that $M$ is a left $A$-module and $N$ is a right $A$-modules. Then $M\o N$ is an $A$-bimodule with the following structures:
\begin{eqnarray}
&&a*(m\o n)=a*m\o n=(h^{a, m}\o 1)\o_H (c_{a, m}\o n)
\end{eqnarray}
and
\begin{eqnarray}
&&(m\o n)*a=m\o n*a=(1\o l^{n, a})\o_H(m\o e_{n, a}),
\end{eqnarray}
where $a* m=(h^{a, m}\o 1)\o_H c_{a, m}, n* a=(1\o l^{n, a})\o_H e_{n, a}$, for all $a\in A, m\in M$ and $n\in N$.

{\bf Proof.} For all $h\in H, m\in M$ and $n\in N$, it is easy to prove that $M\o N$ is a left $H$-module with the action $h\c(m\o n)=h_1m\o h_2n$. We first check that $M\o N$ is a left $A$-module. Observe that
\begin{eqnarray*}
fa*gm&=&((f\o g)\o_H 1)(a*m)\\
&=&(fh^{a, m}\o g)\o_H c_{a,m}\\
&=&(fh^{a, m}S(g_1)\o 1)\o_H g_2c_{a,m}
\end{eqnarray*}
for all $f, g\in H, a\in A$ and $m\in M$. Then we obtain
\begin{eqnarray*}
fa*g(m\o n)&=&fa*(g_1m\o g_2n)\\
&=&(fh^{a, m}S(g_1)\o 1)\o_H (g_2c_{a,m}\o g_3n)\\
&=&(fh^{a, m}S(g_1)\o 1)\D(g_2)\o_H (c_{a,m}\o n)\\
&=&(fh^{a, m}\o g)\o_H (c_{a,m}\o n)\\
&=&((f\o g)\o_H 1)(a*(m\o n)),
\end{eqnarray*}
which proving the $H$-bilinearity. Now we check the associativity. On the one hand,
\begin{eqnarray*}
a*(b*(m\o n))&=&a*((h^{b, m}\o 1)\o_H (c_{b, m}\o n))\\
&=&(1\o h^{b, m}\o 1)(id\o \D)(h^{a, c_{b, m}}\o 1)\o_H (c_{a, c_{b, m}}\o n)\\
&=&(h^{a, c_{b, m}}\o h^{b, m}\o 1)\o_H (c_{a, c_{b, m}}\o n).
\end{eqnarray*}
On the other hand, suppose that $a*b=f^{a,b}\o g^{a,b}\o_H t_{a,b}$ for all $a, b\in A$, then we have
\begin{eqnarray*}
&&(a* b)*(m\o n)=f^{a, b}h^{t_{a, b}, m}_1\o g^{a,b} h^{t_{a, b}, m}_2\o 1\o_H (c_{t_{a, b}, m}\o n).
\end{eqnarray*}
Since $M$ is a left $A$-module, we have that $a*(b*m)=(a*b)*m$, which is equivalent to
\begin{eqnarray*}
&&h^{a, c_{b, m}}\o h^{b, m}\o 1\o_H c_{a, c_{b, m}}=f^{a, b}h^{t_{a, b}, m}_1\o g^{a, b}h^{t_{a, b}, m}_2\o 1\o_H c_{t_{a, b}, m}.
\end{eqnarray*}
It follows that $a*(b*(m\o n))=(a* b)*(m\o n)$. $M\o N$ is also a right $A$-module by a similar calculation. Since
\begin{eqnarray*}
(a*(m\o n))*b&=&((h^{a, m}\o 1)\o_H (c_{a, m}\o n))*b\\
&=&(h^{a, m}\o 1\o 1)(\D\o id)(1\o l^{n, b})\o_H(c_{a, m}\o e_{n, b})\\
&=&h^{a, m}\o 1\o l^{n, b}\o_H(c_{a, m}\o e_{n, b})
\end{eqnarray*}
and
\begin{eqnarray*}
a*((m\o n)*b)&=&a*((1\o l^{n, b})\o_H (m\o e_{n, b}))\\
&=&(1\o 1\o l^{n, b})(id\o \D)(h^{a, m}\o 1)\o_H(c_{a, m}\o e_{n, b})\\
&=&h^{a, m}\o 1\o l^{n, b}\o_H(c_{a, m}\o e_{n, b}),
\end{eqnarray*}
we have that $(a*(m\o n))*b=a*((m\o n)*b)$. Then the conclusion holds.
$\hfill \square$
\\

{\bf Remark 2.5} (1) Let $(A, *)$ be an associative $H$-pseudoalgebra. Then $A\o A$ is an $A$-bimodule by setting $M=N=A$ in Proposition 2.4. More generally, $A^{\o n}(n> 2)$ is an $A$-bimodule with the following structures:
\begin{eqnarray*}
&&a*(b_1\o\cdots\o b_n)=(h^{a, b_1}\o 1)\o_H (c_{a, b_1}\o b_2\o\cdots\o b_n),\\
&&(b_1\o\cdots\o b_n)*a=(1\o l^{b_n, a})\o_H(b_1\o\cdots\o b_{n-1}\o e_{b_n, a})
\end{eqnarray*}
for all $a, b_i(i=1,2,\cdots, n)\in A$. We write
\begin{eqnarray*}
&&a\bullet(b_1\o\cdots\o b_n)=h^{a, b_1}\o (c_{a, b_1}\o b_2\o\cdots\o b_n),\\
&&(b_1\o\cdots\o b_n)\bullet a=l^{b_n, a}\o(b_1\o\cdots\o b_{n-1}\o e_{b_n, a}).
\end{eqnarray*}

(2) Let $(A, *)$ be an associative $H$-pseudoalgebra. Suppose that $M$ is a left $A$-module and $N$ is a left $H$-module. Similar to the proof of Proposition 2.4, $M\o N$ is a left $A$-module with condition $(2.1)$. In addition, if $M$ is a left $H$-module and $N$ is a left $A$-module, then $M\o N$ is also a left $A$-module with the action
\begin{eqnarray*}
&&a*(m\o n)=m\o a* n=(h^{a, n}\o 1)\o_H(m\o c_{a, n}),
\end{eqnarray*}
where $a* n=(h^{a, n}\o 1)\o_H c_{a, n}$.
\\

Let $V$ and $W$ be two left $H$-modules. Recall that an $H$-pseudolinear map from $V$ to $W$ is a $k$-linear map $\phi: V\lr (H\o H)\o_H W$ such that
$$
\phi(hv)=((1\o h)\o_H 1)\phi(v), \quad \forall~ h\in H, v\in V.
$$
The vector space of all such $\phi$ is denoted by $Chom(V, W)$ and the left action of $H$ on $Chom(V, W)$ is defined by
$$
(h\phi)(v)=((1\o h)\o_H 1)\phi(v).
$$
In the special case $V=W$, we write $Cend(V)=Chom(V, V)$. Everywhere in the paper, unless otherwise specified, we always set $V^{*}=Chom(V, k)$.

Consider the map $\r: Chom(V, W)\o V\lr (H\o H)\o_H W$ given by $\r(\phi\o v)=\phi(v)$. By definition it is $H$-bilinear, therefore it is a polylinear map in $\mathcal{M^{*}}(H)$ (see \cite{BDK}). Sometimes, we will use the notation $\phi*v:=\phi(v)$ and consider this as a pseudoproduct or pseudoaction.

 Suppose that $A$ is an associative $H$-pseudoalgebra, $U$ and $V$ are finite $A$-modules. Then the formula
 \begin{eqnarray*}
&&(a*\phi)*u=a*(\phi*u), \quad \forall~ a\in A, u\in U, \phi\in Chom(U, V)
\end{eqnarray*}
provides $Chom(U, V)$ with the structure of a left $A$-module. In particular, when $V$ is the base field $k$, the dual module of $M$ is $M^{*}=Chom(M, k)$, where $k$ is a trival $A$-module with $h\c 1=\v(h)1$ for all $h\in H$.

{\bf Proposition 2.6} Let $M$ and $N$ be two left $A$-modules. Suppose that $M$ is a finitely generated free module (as an $H$-module). Then $M^{*}\o N\simeq Chom(M, N)$ as left $A$-modules, where the correspondence $\phi: M^{*}\o N\lr Chom(M, N)$ is given by
\begin{eqnarray*}
&&\phi(f\o n)*m=(1\o S(g_{f, m}))\o_H n, \forall~ f\in M^{*}, m\in M, n\in N
\end{eqnarray*}
if $f*m=(g_{f, m}\o 1)\o_H 1\in (H\o H)\o_H k$.

{\bf Proof.} Firstly, we check that $\p(f\o n)\in Chom(M, N)$. For all $f\in M^{*}=Chom(M, k)$, we have
\begin{eqnarray*}
f(hm)&=&((1\o h)\o_H 1)f(m)\\
&=&(g_{f, m}\o h)\o_H 1\\
&=&(g_{f, m}S(h_1)\o 1)\o_H h_2\c 1\\
&=&(g_{f, m}S(h)\o 1)\o_H 1.
\end{eqnarray*}
It follows that
\begin{eqnarray*}
\p(f\o n)*(hm)&=&(1\o hS(g_{f, m}))\o_H n\\
&=&((1\o h)\o_H 1)((1\o S(g_{f, m}))\o_H n)\\
&=&((1\o h)\o_H 1)(\p(f\o n)*m).
\end{eqnarray*}
Secondly, we show that $\p$ is a morphism of left $A$-module. By Remark 2.5(2), $M^{*}\o N$ is a left $A$-module via $a*(f\o n)=(h^{a, n}\o 1)\o_H(f\o c_{a, n})$ if $a* n=(h^{a, n}\o 1)\o_H c_{a, n}$.
Since
\begin{eqnarray*}
\p(h(f\o n))*m&=&\p(h_1f\o h_2n)*m\\
&=&(1\o S(h_1g_{f, m}))\o_H h_2n\\
&=&(1\o S(g_{f, m}))(1\o S(h_1))\D(h_2)\o_H n\\
&=&(h\o S(g_{f, m}))\o_H n\\
&=&((h\o 1)\o_H 1)(\p(f\o n)*m)\\
&=&(h\p(f\o n))*m,
\end{eqnarray*}
we have that $\p$ is an $H$-linear map.  For all $m\in M$, we have
\begin{eqnarray*}
(a*\p(f\o n))*m&=&a*(\p(f\o n)*m)\\
&=&a*((1\o S(g_{f, m}))\o_H n)\\
&=&h^{a, n}\o 1\o S(g_{f, m})\o_H c_{a, n}.
\end{eqnarray*}
On the other hand,
\begin{eqnarray*}
\p(a*(f\o n))*m&=&\p((h^{a, n}\o 1)\o_H(f\o c_{a, n}))*m\\
&=&((h^{a, n}\o 1)\o 1)(\D\o id)(1\o S(g_{f, m}))\o_H c_{a, n}\\
&=&£¨h^{a, n}\o 1\o S(g_{f, m}))\o_H c_{a, n}.
\end{eqnarray*}
It follows that $a*\p(f\o n)=\p(a*(f\o n))$.

Finally, it suffices to prove that $\p$ is both injective and surjective. The proofs are similar to Proposition 4.2 in \cite{BL} and we omit the details.
$\hfill \square$

\section{Infinitesimal $H$-pseudobialgebras}\def\theequation{3.\arabic{equation}}
\setcounter{equation} {0}

We start with the following definition.
\\

{\bf Definition 3.1} A {\sl coassociative $H$-coalgebra} $C$ is a left $H$-module, endowed with an $H$-linear map $\D: C\lr C\o C$ ($\D(c)=\sum c_1\o c_2$) satisfying the coassociativity
\begin{eqnarray}
(id\o\D)\circ\D=(\D\o id)\circ\D,
\end{eqnarray}
that is, for all $c\in C$,
\begin{eqnarray*}
&&\sum c_1\o c_{21}\o c_{22}=\sum c_{11}\o c_{12}\o c_2.
\end{eqnarray*}
where $C\o C$ is a left $H$-module via $h\c(c\o d)=h_1c\o h_2d$, for all $h\in H$ and $c, d\in C$.

The coassociative $H$-coalgebra $C$ is cocommutative if it satisfies $\D=\D^{op}$, where $\D^{op}(c)=\sum c_2\o c_1$ for all $c\in C$. For convenience, we omit the summation symbols.
\\

{\bf Remark 3.2} This is nothing but the standard definition of coassociative coalgebra when $H=k$.

Let $L$ be a finite free $H$-module with a basis $\{a_{i}\}_{i=1}^{n}$. The dual basis to $\{a_{i}\}_{i=1}^{n}$ in $L^{*}=Chom(L, k)$ is defined as the set $\{a^{j}\}_{j=1}^{n}$, where each $a^{j}\in L^{*}$ is given by
$$
a^{j}*a_{i}=(1\o 1)\o_H \delta_{ij}.
$$
Obviously, $\{a^{j}\}_{j=1}^{n}$ is a linearly independent set that $H$-generates $L^{*}$.
\\

{\bf Theorem 3.3}(1) Let $(A=\bigoplus_{i=1}^{N}Ha_{i}, *)$ be a finite free associative $H$-pseudoalgebra with the following pseudoproduct:
$$
a_{i}*a_{j}=\sum_{k=1}^{N}(f_{k}^{ij}\o g_{k}^{ij})\o_H a_{k}.
$$
Let $A^{*}=Chom(A, k)=\bigoplus_{i=1}^{N}Ha^{i}$ be the dual of $A$, where $\{a^{i}\}$ is the dual basis corresponding to $\{a_{i}\}$. Define $\D: A^{*}\lr A^{*}\o A^{*}$ as follows:
$$
\D(a^{k})=\sum_{i, j}S(f_{k}^{ij})a^{i}\o S(g_{k}^{ij})a^{j}
$$
and extend it $H$-linearly, i.e., $\D(ha^{k})=h\D(a^{k})$. Then $(A^{*}, \D)$ is a coassociative $H$-coalgebra.

(2) Conversely, suppose $(C, \D)$ is a finite coassociative $H$-coalgebra, then the left $H$-module $C^{*}=Chom(C, k)$ is an associative $H$-conformal algebra with the $x$-product defined by
$$
(f\circ_xg)\circ_y(c)=f\circ_{x_2}(c_1)g\circ_{yS(x_1)}(c_2)
$$
for all $f, g\in C^{*}, c\in C$ and $x, y\in X$.

{\bf Proof.} (1) For all $a^{s}\in L^{*}$, on the one hand,
\begin{eqnarray*}
(id\o\D)\D(a^{s})&=&(id\o\D)(\sum_{i, k}S(f_s^{ik})a^{i}\o S(g_s^{ik})a^{k})\\
&=&\sum_{i, j, l, k}S(f_s^{ik})a^{i}\o\D(S(g_s^{ik}))(S(f_k^{jl}a^{j})\o S(g_k^{jl})a^{l})\\
&=&\sum_{i, j, l, k}(S\o S\o S)(f_s^{ik}\o f_k^{jl}(g_s^{ik})_1\o g_k^{jl}(g_s^{ik})_2)(a^{i}\o a^{j}\o a^{l}).
\end{eqnarray*}
On the other hand,
\begin{eqnarray*}
(\D\o id)\D(a^{s})&=&(\D\o id)(\sum_{k, l}S(f_s^{kl})a^{k}\o S(g_s^{kl})a^{l})\\
&=&\sum_{i, j, l, k}\D(S(f_s^{kl}))(S(f_k^{ij})a^{i}\o S(g_k^{ij})a^{j})\o S(g_s^{kl})a^{l}\\
&=&\sum_{i, j, l, k}(S\o S\o S)(f_k^{ij}(f_s^{kl})_1\o g_k^{ij}(f_s^{kl})_2\o g_s^{kl})(a^{i}\o a^{j}\o a^{l}).
\end{eqnarray*}
Using the associativity of $A$, we have $(a_i*a_j)*a_l=a_i*(a_j*a_l)$, which is equivalent to
\begin{eqnarray}
&&\sum_{k,s}f_k^{ij}(f_s^{kl})_1\o g_k^{ij}(f_s^{kl})_2\o g_s^{kl}\o_H a_s=\sum_{k,s}f_s^{ik}\o f_k^{jl}(g_s^{ik})_1\o g_k^{jl}(g_s^{ik})_2\o_H a_s.
\end{eqnarray}
From $(3.2)$ it follows that $(id\o\D)\D(a^{s})=(\D\o id)\D(a^{s})$. Hence $(A^{*}, \D)$ is a coassociative $H$-coalgebra.

(2) We only prove the associativity of $C^{*}$, the remaining part is similar to Theorem 4.5 in \cite{BL} and we omit the details. For all $f, g, l\in C^{*}$ and $c\in C$, we have
\begin{eqnarray}
\nonumber(f\circ_x(g\circ_yl))\circ_z(c)&=&f\circ_{x_2}(c_1)(g\circ_yl)\circ_{zS(x_1)}(c_2)\\
&=&f\circ_{x_2}(c_1)g\circ_{y_2}(c_{21})l\circ_{(zS(x_1))S(y_1)}(c_{22})
\end{eqnarray}
and
\begin{eqnarray}
\nonumber((f\circ_{x_2}g)\circ_{yx_1}l)\circ_z(c)&=&(f\circ_{x_2}g)\circ_{(yx_1)_2}(c_1)l\circ_{zS((yx_1)_1)}(c_2)\\
&=&f\circ_{x_{22}}(c_{11})g\circ_{(yx_1)_2S(x_{21})}(c_{12})l\circ_{zS((yx_1)_1)}(c_2)
\end{eqnarray}

By using the coassociativity of $\D$ and comparing $(3.3)$ and $(3.4)$, we only need to prove that
\begin{eqnarray}
&&x_{22}\o (yx_1)_2S(x_{21})\o zS((yx_1)_1)=x_2\o y_2\o(zS(x_1))S(y_1)
\end{eqnarray}
for all $x, y, z\in X$. Since $\D(xy)=x_1y_1\o x_2y_2$, we have
\begin{eqnarray*}
&&x_{22}\o(yx_1)_2S(x_{21})\o zS((yx_1)_1)\\
&=&(1\o y_2\o zS(y_1))(x_{22}\o x_{12}S(x_{21})\o S(x_{11}))\\
&=&x_2\o y_2\o(zS(x_1))S(y_1),
\end{eqnarray*}
finishing the proof.
$\hfill \square$
\\
Now we introduce the notion of infinitesimal $H$-pseudobialgebra, which is an associative analogy of Lie $H$-pseudobialgebra.

{\bf Definition 3.4} An {\sl infinitesimal $H$-pseudobialgebra} is a triple $(A, *, \D)$ such that $(A, *)$ is an associative $H$-pseudoalgebra, $(A, \D)$ is a coassociative $H$-coalgebra and they satisfy the compatible condition
\begin{eqnarray}
&&\D(a* b)=a*\D(b)+\D(a)*b, \quad \forall~ a, b\in A,
\end{eqnarray}
where
\begin{eqnarray*}
&&a*\D(b)=\sum_{i}(f_i\o g_i)\o_H (e_i\o b_2),\\
&&\D(a)*b=\sum_{j}(k_j\o l_j)\o_H(a_1\o t_j),
\end{eqnarray*}
if $a*b_1=\sum_{i}(f_i\o g_i)\o_H e_i, a_2*b=\sum_{j}(k_j\o l_j)\o_H t_j$.
\\

{\bf Remark 3.5} In particular, for the one dimensional Hopf algebra $H=k$, an infinitesimal $H$-pseudobialgebra is an ordinary infinitesimal bialgebra (\cite{A4}) over the field $k$.
\\

{\bf Example 3.6} Let $A=H\{e_1, e_2\}$ be a free associative $H$-pseudoalgebra with the pseudoproduct given by $e_1*e_2=e_2*e_1=e_2*e_2=0, e_1*e_1=(f\o g)\o_H e_2$, $\forall~f, g\in H$. Define $\D: A\lr A\o A$
as follows:
$$
\D(e_1)=e_1\o e_2,~\D(e_2)=e_2\o e_2,
$$
and extend it $H$-linearity, i.e., $\D(he_i)=h\D(e_i)$ for $i=1, 2$. Then $(A,\D)$ is a coassociative $H$-coalgebra. Furthermore, $(A,*,\D)$ is an infinitesimal $H$-pseudobialgebra. Suppose that $r=e_i\o e_j, i,j=1, 2$ except for $i=j=1$, then $(A,*,\D_{r})$ is a coboundary infinitesimal $H$-pseudobialgebra, which will be defined in the next section.
\\

{\bf Proposition 3.7} Let $H'$ be a Hopf subalgebra of $H$ and $(A, *, \D)$ an infinitesimal $H'$-pseudobialgebra. Then $(Cur(A)=H\o_{H'} A, \tilde{*}, \delta)$ is an infinitesimal $H$-pseudobialgebra with the following structures ($\forall~ f, g\in H, a, b\in A$):
\begin{eqnarray*}
&&(f\o_{H'} a)\tilde{*}(g\o_{H'} b)=\sum_{i}(ff_i\o gg_i)\o_H(1\o_{H'} e_i),\\
&&\delta(f\o_{H'} a)=(f_1\o_{H'} a_1)\o (f_2\o_{H'} a_2),
\end{eqnarray*}
where $a*b=\sum_{i}(f_i\o g_i)\o_{H'} e_i$.

{\bf Proof.} By \cite{BDK}, $(Cur(A)=H\o_{H'} A, \tilde{*})$ is an associative $H$-pseudoalgebra. For all $f, g\in H$ and $a, b\in A$, we have
\begin{eqnarray*}
&&\delta(f\c(g\o_{H'} b))=\delta(fg\o_{H'} b)=(f_1g_1\o_{H'} b_1)\o(f_2g_2\o_{H'} b_2)=f\c\delta(g\o_{H'} b)
\end{eqnarray*}
and
\begin{eqnarray*}
(id\o\d)\d(f\o_{H'} a)&=&(id\o\d)((f_1\o_{H'} a_1)\o(f_2\o_{H'} a_2))\\
&=&((f_1\o_{H'} a_1))\o((f_2\o_{H'} a_{21}))\o((f_3\o_{H'} a_{22}))\\
&=&((f_1\o_{H'} a_{11}))\o((f_2\o_{H'} a_{12}))\o((f_3\o_{H'} a_{2}))\\
&=&(\d\o id)\d(f\o_{H'} a).
\end{eqnarray*}
It follows that $(Cur(A), \d)$ is a coassociative $H$-coalgebra. In what follows, we verify the compatible condition (3.6). Suppose
$$
a*b=\sum_{i}(f_i\o g_i)\o_{H'} e_i,~a*b_1=\sum_{i}(l_i\o k_i)\o_{H'}t_i,~a_2*b=\sum_{i}(m_i\o n_i)\o_{H'}r_i.
$$
We compute:
\begin{eqnarray*}
\d((f\o_{H'} a)\tilde{*}(g\o_{H'} b))&=&\d(\sum_{i}(ff_i\o gg_i)\o_H(1\o_{H'}e_i))\\
&=&\sum_{i}(ff_i\o gg_i)\o_H((1\o_{H'} e_{i1})\o(1\o_{H'} e_{i2})),
\end{eqnarray*}
\begin{eqnarray*}
(f\o_{H'} a)\tilde{*}\d(g\o_{H'} b)&=&(f\o_{H'} a)\tilde{*}(g_1\o_{H'} b_1)\o (g_2\o_{H'} b_2)\\
&=&\sum_{i}(fl_i\o g_1k_i)\o_H((1\o_{H'} t_i)\o(g_2\o_{H'} b_2))\\
&=&\sum_{i}(fl_i\o gk_i)\o_H((1\o_{H'} t_i)\o(1\o_{H'} b_2))
\end{eqnarray*}
and
\begin{eqnarray*}
\d(f\o_{H'} a)\tilde{*}(g\o_{H'} b)&=&(f_1\o_{H'} a_1)\o(f_2\o_{H'} a_2)\tilde{*}(g\o_{H'} b)\\
&=&\sum_{i}(f_2m_i\o gn_i)\o_H((f_1\o_{H'} a_1)\o(1\o_{H'} r_i))\\
&=&\sum_{i}(fm_i\o gn_i)\o_H((1\o_{H'} a_1)\o(1\o_{H'} r_i)).
\end{eqnarray*}
Since $(A, *, \D)$ is an infinitesimal $H'$-pseudobialgebra, we have
$$
\D(a* b)=\D(a)* b+a*\D(b),
$$
that is,
$$
\sum_{i}(f_i\o g_i)\o_{H'}(e_{i1}\o e_{i2})=\sum_{i}(l_i\o k_i)\o_{H'}(t_i\o b_2)+\sum_{i}(m_i\o n_i)\o_{H'}(a_1\o r_i),
$$
which implies that $\d((f\o_{H'} a)\tilde{*}(g\o_{H'} b))=(f\o_{H'} a)\tilde{*}\d(g\o_{H'} b)+\d(f\o_{H'} a)\tilde{*}(g\o_{H'} b)$. This completes the proof.
$\hfill \square$
\\

{\bf Remark 3.8} More generally, let $\phi: H'\lr H$ be a homomorphism of Hopf algebras, and $(A, *, \D)$ be an infinitesimal $H'$-pseudobialgebra. Then $(H\o_{H'}A, \tilde{*}, \d)$ is an infinitesimal $H$-pseudobialgebra with the strucutres:
\begin{eqnarray*}
&&(f\o_{H'} a)\tilde{*}(g\o_{H'} b)=\sum_{i}(f\phi(f_i)\o g\phi(g_i))\o_H (1\o_{H'}e_i),\\
&&\d(f\o_{H'} a)=(f_1\o_{H'} a_1)\o(f_2\o_{H'} a_2),
\end{eqnarray*}
for all $f, g\in H$ and $a, b\in A$.
\\

{\bf Corollary 3.9} Let $(A, \mu, \D)$ be an infinitesimal bialgebra. Then $(Cur(A)=H\o A, \tilde{*}, \delta)$ is an infinitesimal $H$-pseudobialgebra with the following structures:
\begin{eqnarray*}
&&(f\o a)\tilde{*}(g\o b)=(f\o g)\o_H(1\o ab),\\
&&\delta(f\o a)=(f_1\o a_1)\o (f_2\o a_2),
\end{eqnarray*}
for all $f, g\in H$ and $a, b\in A$.

{\bf Proof.} It can be obtained directly by taking $H'=k$ in Proposition 3.7.
$\hfill \square$

\section{Coboundary infinitesimal $H$-pseudobialgebras}\def\theequation{4.\arabic{equation}}
\setcounter{equation} {0}

In this section, we study an important subclass of infinitesimal $H$-pseudobialgebra, for which the coalgebra structure comes from a 1-coboundary in associative $H$-pseudoalgebra cohomology. First we introduce the cohomology of associative $H$-pseudoalgebras.

Let $A$ be an associative $H$-pseudoalgebra and $M$ an $A$-bimodule. Define $C^{n}(A, M)$$(n\geq 1)$, consisting of all cochains
$$
\gamma\in Hom_{H^{\o n}}(A^{\o n}, H^{\o n}\o_H M).
$$
Explicitly, $\gamma$ has the following defining property: $H$-polylinearity,
$$
\g(h_1a_1\o\cdots\o h_na_n)=((h_1\o\cdots\o h_n)\o_H 1)\g(a_1\o\cdots\o a_n)
$$
where $h_i\in H$ and $a_i\in A$ for $i=1, 2,\cdots, n$.

For $n=0$, we put $C^{0}(A, M)=k\o_H M\simeq M/H_{+}M$, where $H_{+}=\{h\in H|\v(h)=0\}$ is the augmentation ideal of $H$. The differential $d: C^{0}(A, M)\lr C^{1}(A, M)=Hom_H(A, M)$ is given by
\begin{eqnarray*}
(d(1\o_H m))(a)=\sum_{i}(1\o\v)(h_i)m_i-\sum_{j}(\v\o 1)(f_j)n_j, \forall~ a\in A, m\in M,
\end{eqnarray*}
if $a*m=\sum_{i}h_i\o_H m_i\in H^{\o 2}\o_H M$ and $m* a=\sum_{j}f_j\o_H n_j\in H^{\o 2}\o_H M$.

For $n\geq 1$, the differential $d: C^{n}(A, M)\lr C^{n+1}(A, M)$ is given by
\begin{eqnarray}
\nonumber&&d\g(a_1\o\cdots\o a_{n+1})=\r_l(a_1\o\g(a_2\o\cdots\o a_{n+1}))\\
\nonumber&&+\sum_{i=1}^{n}(-1)^{i}\g(a_1\o\cdots\o a_{i-1}\o\mu(a_i\o a_{i+1})\o a_{i+2}\o\cdots\o a_{n+1})\\
&&+(-1)^{n+1}\r_r(\g(a_1\o\cdots\o a_n)\o a_{n+1}).
\end{eqnarray}

We also use the following convention in the above equation. If $a*m=\sum_{i}h_i\o_H m_i\in H^{\o 2}\o_H M$, $m*a=\sum_{j}f_j\o_H n_j\in H^{\o 2}\o_H M$ for all $a\in A$ and $m\in M$, then for any $g\in H^{\o n}$, we set
\begin{eqnarray*}
&&a*(g\o m)=\sum_{i}(1\o g\D^{(n-1)})(h_i)\o_H m_i\in H^{\o n+1}\o_H M\\
\hbox{and}\\
&&(g\o m)*a=\sum_{j}(g\D^{(n-1)}\o 1)(f_j)\o_H n_j\in H^{\o n+1}\o_H M,
\end{eqnarray*}
where $\D^{(n-1)}=(id\o\cdots\o id\o\D)\cdots(id\o\D)\D: H\lr H^{\o n}$ is the iterated comultiplication ($\D^{(0)}=id$). Note that equation $(4.1)$ holds also for $n=0$ if we define $\D^{(-1)}=\v$.

Equation $(4.1)$ is illustrated in Figure 1.

\begin{figure*}[!htb]
  \centering
  \includegraphics[width=\textwidth]{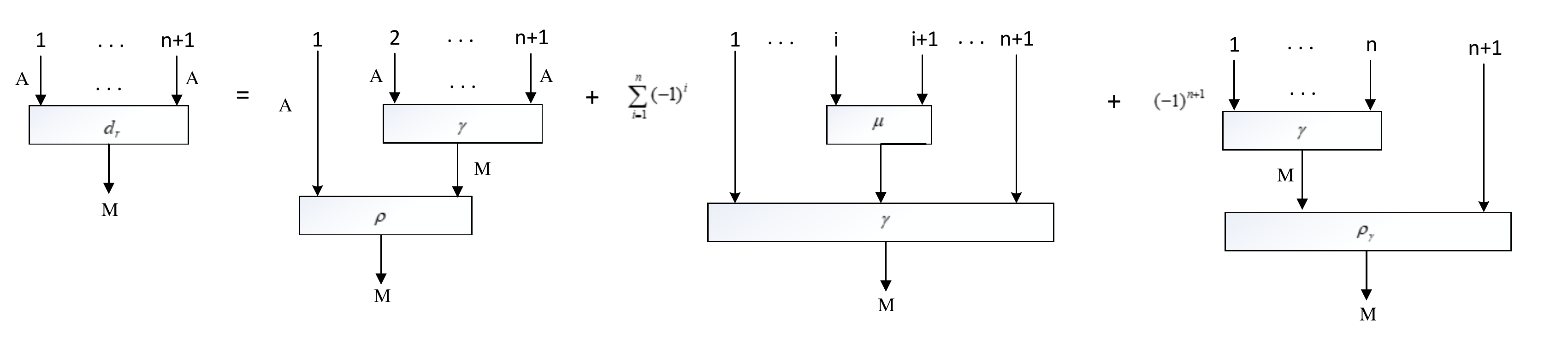}
  \caption{The definition of differential}
  \label{fig1}
\end{figure*}

One can verify that $d^{2}=0$ by using the same argument as in the usual associative algebra case. $\gamma\in C^{n}(A, M)$ is called an $n$-cocycle if $d \gamma=0$. If $\gamma=d \alpha$ for some $\a\in C^{n-1}(A, M)$, then $\gamma$ is called an $n$-coboundary. The cohomology of the resulting complex $C^{\bullet}(A, M)$ is called the reduced cohomology of $A$ with coefficient in $M$ and is denoted by $H^{\bullet}(A, M)$.
\\

{\bf Remark 4.1} For $n=1$, $\gamma\in C^{1}(A, M)$ is a 1-cocycle if $\gamma(a_1*a_2)=a_1*\gamma(a_2)+\gamma(a_1)*a_2$. Now if $A$ is an infinitesimal $H$-pseudobialgebra with comultiplication $\D$, then the compatible condition $(3.6)$ is indeed the condition that $\D$ is a 1-cocycle of $A$ with coefficient in $A\o A$ in the reduced complex.

Among all the 1-cocycles of $A$ with values in $A\o A$, we have 1-coboundary $\D_r$ that comes from the differential of an element $r\in A\o A$, that is $\D_r(a)=(d(1\o_H r))(a)$ for all $a\in A$. More precisely, suppose $r=\sum_{i}u_i\o v_i\in A\o A$, then
\begin{eqnarray}
\D_r(a)=d(1\o_H r)(a)=h^{a, u_i}\c(c_{a, u_i}\o v_i)-l^{v_i, a}\c(u_i\o e_{v_i, a}),
\end{eqnarray}
where $a*u_i=h^{a, u_i}\o 1\o_H c_{a, u_i}, v_i*a=1\o l^{v_i, a}\o_H e_{v_i, a}$.

Now, we introduce the subclass of infinitesimal $H$-pseudobialgebras.
\\

{\bf Definition 4.2} A {\sl coboundary infinitesimal $H$-pseudobialgebra} is a quadruple $(A, *, \D_r, r)$, with $r\in A\o A$, such that $(A, *, \D_r)$ is an infinitesimal $H$-pseudobialgebra.
\\

{\bf Example 4.3} Let $A=H\{e_1, e_2\}$ be a free associative $H$-pseudoalgebra with pseudoproduct given by
$$
e_1*e_1=e_1*e_2=0,~~e_2*e_1=1\o 1\o_H e_1,~~e_2*e_2=1\o 1\o_H e_2.
$$
 Define$r=e_2\o e_1-e_1\o e_2\in A\o A$. Then by straightforward computations, $\delta_r(e_1)=e_1\o e_1, \delta_r(e_2)=e_2\o e_1$, and $(A, *, r)$ is a coboundary infinitesimal $H$-pseudobialgebra.
\\

Let $(A, *)$ be an associative $H$-pseudoalgebra and $r=\sum_{i} u_i\o v_i\in A\o A$. We define the associative pseudo-Yang-Baxter equation ($pseudo$-AYBE) on $A$ as
\begin{eqnarray}
\nonumber A(r)&=&\mu_{-1}^{4}(\{u_i, u_j\}\o v_j\o v_i)-\mu_{2}^{3,4}(u_i\o \{v_i, u_j\}\o v_j)\\
&&+\mu_{3}^{1, 4}(u_i\o u_j\o \{v_j, v_i\})
\end{eqnarray}
in $A\o A\o A$, where $\mu_{-k}^{l}$ means that the element of $H$ that appears in its argument in the $k$-th place acts via the antipode on the element of $A$ located in the $l$-th entry, $\mu_{k}^{r,s}$ means that the element of $H$ located in the $k$-th place in its argument acts on the elements of $A\o A$ formed by the elements in the $r$-th and $s$-th places. For example, $\mu_{-1}^{4}(h\o a\o b\o c)=a\o b\o S(h)c$, $\mu_{3}^{1,4}(a\o b\o h\o c)=h_1a\o b\o h_2c$ for all $h\in H$ and $a, b, c\in A$.

We say that $r$ is a solution of the $pseudo$-AYBE if $A(r)=0$. Moreover, the $pseudo$-AYBE is exactly the usual AYBE when $H=k$.

In a coboundary infinitesimal $H$-pseudobialgebra, the comultiplication $\D_r$ is determined by $r$. If $(A, *)$ is an associative $H$-pseudoalgebra and $r\in A\o A$. We discuss what conditions $(A, *, \D_r, r)$ is a coboundary infinitesimal $H$-pseudobialgebra in the following.

{\bf Theorem 4.4} Let $(A, *)$ be an associative $H$-pseudoalgebra and $r=\sum_{i}u_i\o v_i\in A\o A$. Then $(A, *, \D_r, r)$ is a coboundary infinitesimal $H$-pseudobialgebra if and only if
\begin{eqnarray}
\mu_{3}(a\bullet A(r)-A(r)\bullet a)=0,
\end{eqnarray}
where $\mu_3(h\o a\o b\o c)=((\D\o id)\D(h))(a\o b\o c)$ for all $h\in H$ and $a, b, c\in A$.

{\bf Proof.} By the definition of $\D_r$, we have that $\D_r$ satisfies the compatible condition $(3.6)$ regardless of whether condition $(4.4)$ holds or not. So we only need to prove that $\D_r$ is coassociativity if and only if $\mu_{3}(a\bullet A(r)-A(r)\bullet a)=0$.

Suppose $a* u_i=h^{a, u_i}\o 1\o_H c_{a, u_i}, v_i*a=1\o l^{v_i, a}\o_H e_{v_i, a}$.
Then $\D_r(a)=h^{a, u_i}\c(c_{a, u_i}\o v_i)-l^{v_i, a}\c(u_i\o e_{v_i, a})$. On the one hand, we compute:
\begin{eqnarray}
\nonumber&&(\D_r\o id)\D_r(a)-(id\o\D_r)\D_r(a)\\
\nonumber&=&(\D_r\o id)(h^{a, u_i}\c(c_{a, u_i}\o v_i)-l^{v_i, a}\c(u_i\o e_{v_i, a}))\\
\nonumber&&-(id\o\D_r)(h^{a, u_i}\c(c_{a, u_i}\o v_i)-l^{v_i, a}\c(u_i\o e_{v_i, a}))\\
\nonumber&=&h^{a, u_i}_1\c\D_r(c_{a, u_i})\o h^{a, u_i}_2 v_i-l^{v_i, a}_1\c\D_r(u_i)\o l^{v_i, a}_2 e_{v_i, a}\\
\nonumber&&-h^{a, u_i}_1 c_{a, u_i}\o h^{a, u_i}_2\c \D_r(v_i)+l^{v_i, a}_1 u_i\o l^{v_i, a}_2\c \D_r(e_{v_i, a})\\
\nonumber&=&h^{a, u_i}_1\c(h^{c_{a, u_i}, u_j}\c(c_{c_{a, u_i}, u_j}\o v_j)-l^{v_j, c_{a, u_i}}\c(u_j\o e_{v_j, c_{a, u_i}}))\o h^{a, u_i}_2 v_i\\
\nonumber&&-l^{v_i, a}_1\c(h^{u_i, u_j}\c(c_{u_i, u_j}\o v_j)-l^{v_j, u_i}\c(u_j\o e_{v_j, u_i}))\o l^{v_i, a}_2 e_{v_i, a}\\
\nonumber&&-h^{a, u_i}_1 c_{a, u_i}\o h^{a, u_i}_2\c (h^{v_i, u_j}\c(c_{v_i, u_j}\o v_j)-l^{v_j, v_i}\c(u_j\o e_{v_j, v_i}))\\
\nonumber&&+l^{v_i, a}_1 u_i\o l^{v_i, a}_2\c (h^{e_{v_i, a}, u_j}\c(c_{e_{v_i, a}, u_j}\o v_j)-l^{v_j, e_{v_i, a}}\c(u_j\o e_{v_j, e_{v_i, a}}))\\
&=&h^{a, u_i}\c(h^{c_{a, u_i}, u_j}\c(c_{c_{a, u_i}, u_j}\o v_j)\o v_i)\\
&&-h^{a, u_i}\c(l^{v_j, c_{a, u_i}}\c(u_{j}\o e_{v_{j}, c_{a, u_i}})\o v_i)\\
&&-l^{v_i, a}\c(h^{u_i, u_j}\c(c_{u_i, u_j}\o v_{j})\o e_{v_i, a})\\
&&+l^{v_i, a}\c(l^{v_j, u_i}\c(u_j\o e_{v_j, u_i})\o e_{v_i, a})\\
&&-h^{a, u_i}\c(c_{a, u_i}\o h^{v_i, u_j}\c(c_{v_i, u_j}\o v_j))\\
&&+h^{a, u_i}\c(c_{a, u_i}\o l^{v_j, v_i}\c(u_j\o e_{v_j, v_i}))\\
&&+l^{v_i, a}\c(u_i\o h^{e_{v_i, a}, u_j}\c(c_{e_{v_i, a}, u_j}\o v_j))\\
&&-l^{v_i, a}\c(u_i\o l^{v_j, e_{v_i, a}}\c(u_j\o e_{v_j, e_{v_i, a}})).
\end{eqnarray}
On the other hand, we have
\begin{eqnarray}
\nonumber&&\mu_3(a\bullet A(r)-A(r)\bullet a)\\
&=&\mu_3(a\bullet(\mu_{-1}^{4}(\{u_i, u_j\}\o v_j\o v_i)))\\
&&-\mu_3(a\bullet(\mu_{2}^{3, 4}(u_i\o \{v_i, u_j\}\o v_j)))\\
&&+\mu_3(a\bullet(\mu_{3}^{1,4}(u_i\o u_j\o \{v_j, v_i\})))\\
&&-\mu_3((\mu_{-1}^{4}(\{u_i, u_j\}\o v_j\o v_i))\bullet a)\\
&&+\mu_3((\mu_{2}^{3,4}(u_i\o \{v_i, u_j\}\o v_j))\bullet a)\\
&&-\mu_3((\mu_{3}^{1,4}(u_i\o u_j\o \{v_j, v_i\}))\bullet a).
\end{eqnarray}
Next, we only need to verify that
 \begin{eqnarray}
 (\D_r\o id)\D_r(a)-(id\o\D_r)\D_r(a)-\mu_3(a\bullet A(r)-A(r)\bullet a)=0.
 \end{eqnarray}
 Using the property of Fourier transform, we obtain
\begin{eqnarray*}
&&\{\{a, u_i\}, u_j\}=(\mathcal{F}\o id)\{a, \{u_i, u_j\}\}\\
&=&(\mathcal{F}\o id)(h^{a, c_{u_i, u_j}}\o h^{u_i, u_j}\o_H c_{a, c_{u_i, u_j}})\\
&=&h^{a, c_{u_i, u_j}}S(h^{u_i, u_j}_1)\o h^{u_i, u_j}_2\o_H c_{a, c_{u_i, u_j}}.
\end{eqnarray*}
So we have
\begin{eqnarray*}
&&h^{a, u_i}\c(h^{c_{a, u_i}, u_j}\c(c_{c_{a, u_i}, u_j}\o v_j)\o v_i)\\
&=&\mu_3(\mu_{2}^{3, 4}(\{\{a, u_i\}, u_j\}\o v_j\o v_i))\\
&=&\mu_3(\mu_{2}^{3, 4}((\mathcal{F}\o id)\{a, \{u_i, u_j\}\}\o v_j\o v_i))\\
&=&(h^{a, c_{u_i, u_j}}S(h^{u_i, u_j}_1))\c(h^{u_i, u_j}_2\c(c_{a, c_{u_i, u_j}}\o v_j)\o v_i)\\
&=&h^{a, c_{u_i, u_j}}\c(c_{a, c_{u_i, u_j}}\o v_j\o S(h^{u_i, u_j})v_i)\\
&=&\mu_3(a\bullet(\mu_{-1}^{4}(\{u_i, u_j\}\o v_j\o v_i))).
\end{eqnarray*}
Hence $(4.5)-(4.13)=0$. Since
\begin{eqnarray*}
S(h^{u_i, u_j})v_i*a&=&(S(h^{u_i, u_j})\o 1)(1\o l^{v_i, a})\o_H e_{v_i, a}\\
&=&S(h^{u_i, u_j})\o l^{v_i, a}\o_H e_{v_i, a}\\
&=&1\o l^{v_i, a}h^{u_i, u_j}_1\o_H S(h^{u_i, u_j}_2)e_{v_i, a},
\end{eqnarray*}
we have
\begin{eqnarray*}
&&l^{v_i, a}\c(h^{u_i, u_j}\c(c_{u_i, u_j}\o v_{j})\o e_{v_i, a})\\
&=&(l^{v_i, a}h^{u_i, u_j}_1)\c(c_{u_i, u_j}\o v_{j}\o S(h^{u_i, u_j}_2)e_{v_i, a})\\
&=&\mu_3((c_{u_i, u_j}\o v_{j}\o S(h^{u_i, u_j})v_i)\bullet a)\\
&=&\mu_3((\mu_{-1}^{4}(\{u_i, u_j\}\o v_j\o v_i))\bullet a).
\end{eqnarray*}
It follows that $(4.7)-(4.16)=0$.
Similarly, we get
\begin{eqnarray*}
&&l^{v_i, a}\c(l^{v_j, u_i}\c(u_j\o e_{v_j, u_i})\o e_{v_i, a})\\
&=&(l^{v_i, a}l^{v_j, u_i}_1)\c(u_j\o e_{v_j, u_i}\o S(l^{v_j, u_i}_2)e_{v_i, a})\\
&=&\mu_3((u_j\o e_{v_j, u_i}\o S(l^{v_j, u_i})v_i)\bullet a)\\
&=&\mu_3(\mu_{2}^{3,4}(u_j\o S(l^{v_j, u_i}_1)\o l^{v_j, u_i}_2e_{v_j, u_i}\o v_i)\bullet a)\\
&=&\mu_3((\mu_{2}^{3,4}(u_i\o \{v_i, u_j\}\o v_j))\bullet a)
\end{eqnarray*}
and
\begin{eqnarray*}
&&h^{a, u_i}\c(c_{a, u_i}\o l^{v_j, v_i}\c(u_j\o e_{v_j, v_i}))\\
&=&(h^{a, u_i}l^{v_j, v_i}_1)\c(S(l^{v_j, v_i}_2)c_{a, u_i}\o u_j\o e_{v_j, v_i})\\
&=&\mu_3(a\bullet(S(l^{v_j, v_i})u_i\o u_j\o e_{v_j, v_i}))\\
&=&\mu_3(a\bullet \mu_{3}^{1,4}(u_i\o u_j\o S(l^{v_j, v_i}_1)\o l^{v_j, v_i}_2e_{v_j, v_i}))\\
&=&\mu_3(a\bullet(\mu_{3}^{1,4}(u_i\o u_j\o \{v_j, v_i\}))).
\end{eqnarray*}
Thus $(4.8)-(4.17)=0$ and $(4.10)-(4.15)=0$. Clearly,
\begin{eqnarray*}
&&h^{a, u_i}\c(c_{a, u_i}\o h^{v_i, u_j}\c(c_{v_i, u_j}\o v_j))\\
&=&\mu_3(a\bullet(u_i\o h^{v_i, u_j}\c(c_{v_i, u_j}\o v_j)))\\
&=&\mu_3(a\bullet(\mu_{2}^{3, 4}(u_i\o \{v_i, u_j\}\o v_j))).
\end{eqnarray*}
So we have $(4.9)-(4.14)=0$. Using the associativity of $A$, we have $(v_i* a)*u_j=v_i*(a* u_j)$. Suppose
\begin{eqnarray*}
&&v_i* a=1\o l^{v_i,a}\o_H e_{v_i, a},\\
&&a*u_j=h^{a, u_j}\o 1\o_H c_{a, u_j}=1\o S(h^{a, u_j}_1)\o_H h^{a, u_j}_2c_{a, u_j},\\
&&e_{v_i, a}*u_j=h^{e_{v_i, a}, u_j}\o 1\o_H c_{e_{v_i, a}, u_j}=1\o S(h^{e_{v_i, a}, u_j}_1)\o_H h^{e_{v_i, a}, u_j}_2c_{e_{v_i, a}, u_j},\\
&&v_i*c_{a, u_j}=1\o l^{v_i, c_{a, u_j}}\o_H e_{v_i, c_{a, u_j}}.
\end{eqnarray*}
Then we have
\begin{eqnarray*}
&&1\o l^{v_i, a}\o S(h^{e_{v_i, a}, u_j}_1)\o_H h^{e_{v_i, a}, u_j}_2c_{e_{v_i, a}, u_j}\\
&=&1\o h^{a, u_j}l^{v_i, c_{a, u_j}}_1\o l^{v_i, c_{a, u_j}}_2\o_H e_{v_i, c_{a, u_j}}.
\end{eqnarray*}
Applying $\v\o id\o S\o id$ to the above equation, we obtain
$$
g^{v_i, a}\o h^{e_{v_i, a}, u_j}_1\o_H h^{e_{v_i, a}, u_j}_2c_{e_{v_i, a}, u_j}=h^{a, u_j}l^{v_i, c_{a, u_j}}_1\o S(l^{v_i, c_{a, u_j}}_2)\o_H e_{v_i, c_{a, u_j}}.
$$
Thus we get that $(4.6)+(4.11)$ is
\begin{eqnarray*}
&&l^{v_i, a}\c(u_i\o h^{e_{v_i, a}, u_j}\c(c_{e_{v_i, a}, u_j}\o v_j))-h^{a, u_i}\c(l^{v_j, c_{a, u_i}}\c(u_{j}\o e_{v_{j}, c_{a, u_i}})\o v_i)\\
&=&l^{v_i, a}\c(u_i\o h^{e_{v_i, a}, u_j}_2c_{e_{v_i, a}, u_j}\o h^{e_{v_i, a}, u_j}_1v_j)-h^{a, u_i}\c(l^{v_j, c_{a, u_i}}\c(u_{j}\o e_{v_{j}, c_{a, u_i}})\o v_i)\\
&=&(h^{a, u_j}l^{v_i, c_{a, u_j}}_1)\c(u_i\o e_{v_i, c_{a, u_j}}\o S(l^{v_i, c_{a, u_j}}_2)v_j)-h^{a, u_i}\c(l^{v_j, c_{a, u_i}}\c(u_{j}\o e_{v_{j}, c_{a, u_i}})\o v_i)\\
&=&0.
\end{eqnarray*}
Using the associativity of $A$ again, we have $v_j*(v_i* a)=(v_j* v_i)*a$. Suppose
\begin{eqnarray*}
&&v_i* a=1\o l^{v_i,a}\o_H e_{v_i, a},~~~ v_j*e_{v_i, a}=1\o l^{v_j, e_{v_i, a}}\o_H e_{v_j, e_{v_i, a}},\\
&&v_j*v_i=1\o l^{v_j, v_i}\o_H e_{v_j, v_i},~~~ e_{v_j, v_i}*a=1\o l^{e_{v_j, v_i}, a}\o e_{e_{v_j, v_i}, a}.
\end{eqnarray*}
Then we obtain
\begin{eqnarray*}
&&1\o l^{v_j, e_{v_i, a}}_1\o l^{v_i, a}l^{v_j, e_{v_i, a}}_2\o_H e_{v_j, e_{v_i, a}}\\
&=&1\o l^{v_j, v_i}\o l^{e_{v_j, v_i}, a}\o_H e_{e_{v_j, v_i}, a}.
\end{eqnarray*}
It follows that
\begin{eqnarray*}
&&l^{v_i, a}\c(u_i\o l^{v_j, e_{v_i, a}}\c(u_j\o e_{v_j, e_{v_i, a}}))\\
&=&(l^{v_i, a}l^{v_j, e_{v_i, a}}_2)\c(S(l^{v_j, e_{v_i, a}}_1)u_i\o u_j\o e_{v_j, e_{v_i, a}})\\
&=&l^{e_{v_j, v_i}, a}\c(S(l^{v_j, v_i})u_i\o u_j\o e_{e_{v_j, v_i}, a})\\
&=&\mu_3((S(l^{v_j, v_i})u_i\o u_j\o e_{v_j, v_i})\bullet a)\\
&=&\mu_3((\mu_{3}^{1,4}(u_i\o u_j\o \{v_j, v_i\}))\bullet a).
\end{eqnarray*}
So $(4.12)-(4.18)=0$. Finally, it is easy to check that we have canceled all the terms of the left-hand side of equation $(4.19)$. This completes the proof.
$\hfill \square$

\section{From (coboundary) infinitesimal $H$-pseudobialgebras to (coboundary) Lie $H$-pseudobialgebras}\def\theequation{5.\arabic{equation}}
\setcounter{equation} {0}

We first recall some definitions about Lie $H$-pseudobialgebras (see \cite{BL}).
\\

{\bf Definition 5.1} A Lie $H$-pseudoalgebra $(L, [*])$ is a left $H$-module $L$ endowed with a map (called the pseudobracket):
\begin{eqnarray*}
[*]: \quad L\o L\lr (H\o H)\o_H L,\quad a\o b\mapsto [a*b]
\end{eqnarray*}
satisfying ($\forall~ a, b, c\in L, f, g\in H$)

{\bf $H$-bilinearity:} \quad $[fa*gb]=(f\o g\o_H 1)[a*b]$.

{\bf Skew-commutativity:} $[a*b]=-(\sigma\o id)[b*a]$.

{\bf Jacobi identity:} $[[a*b]*c]=[a*[b*c]]-((\sigma\o id)\o_H id)[b*[a*c]]$.
\\

{\bf Definition 5.2} A Lie $H$-coalgebra is a left $H$-module, endowed with an $H$-linear map $\delta: L\lr L\o L$ such that
\begin{eqnarray*}
&&\tau\circ\delta=-\delta,\\
&&(id\o\d)\d-\sigma_{12}(id\o\d)\d=(\d\o id)\d,
\end{eqnarray*}
where $\tau: L\o L\lr L\o L$ is the permutation $\tau(a\o b)=b\o a$, $\forall~a, b\in L$.
\\

{\bf Definition 5.3} A Lie $H$-pseudobialgebra is a triple $(L, [*], \d)$ such that $(L, [*])$ is a Lie $H$-pseudoalgebra, $(L, \d)$ is a Lie $H$-coalgebra, and they satisfy the cocycle condition
\begin{eqnarray}
&&\d([a* b])=[a*\d(b)]-(\sigma\o_H id)[b*\d(a)],
\end{eqnarray}
where $[a*\d(b)]=\sum_{i}(f_iS(g_{i1})\o 1)\o_H(g_{i2}a_i\o b_2)+\sum_{j}(k_jS(l_{j1})\o 1)\o_H(b_1\o l_{j2}b_j)$, if
\begin{eqnarray*}
&&[a* b_1]=\sum_{i}(f_i\o g_i)\o_H a_i=\sum_{i}(f_iS(g_{i1})\o 1)\o_H g_{i2}a_i
\end{eqnarray*}
and
\begin{eqnarray*}
&&[a* b_2]=\sum_{j}(k_j\o l_j)\o_H b_j=\sum_{j}(k_jS(l_{j1})\o 1)\o_H l_{j2}b_j.
\end{eqnarray*}
\\

{\bf Definition 5.4} A coboundary Lie $H$-pseudobialgebra $(L, [*], \d, r)$ consists of a Lie $H$-pseudo-\\
bialgebra $(L, [*], \d)$ and an element $r=\sum_{i}u_i\o v_i\in L\o L$ such that
\begin{eqnarray}
\d(a)=\sum_{i}\mu([a, u_i]\o v_i+\sigma_{12}(u_i\o[a, v_i])),
\end{eqnarray}
where $\mu: H\o(L\o L)\lr L\o L$ is given by $\mu(h\o m\o n)=\D(h)(m\o n)$, and $[a, b]$ is the Fourier transform of $[a*b]$.

To state the main result of this section, we need the following notation.

Let $A$ be an infinitesimal $H$-pseudoalgebra.
Define the $\star$-bracket
\begin{eqnarray*}
&&[a, b]_{\star}=a* b-(\sigma\o id)(b* a),
\end{eqnarray*}
for all $a\in A$ and $b=b_1\o\cdots\o b_n\in A^{\o n}$.
\\

{\bf Definition 5.5} Let $(A, *, \D)$ be an infinitesimal $H$-pseudobialgebra. Define the map $\mathfrak{B}: A\o A\lr (H\o H)\o_H (A\o A)$ by
\begin{eqnarray*}
&&\mathfrak{B}(a, b)=[a, \D^{op}(b)]_{\star}+(\sigma\o \tau)([b, \D^{op}(a)]_{\star})
\end{eqnarray*}
 where $\tau(a\o b)=b\o a$, for all $a, b\in A$. We call $\mathfrak{B}$ the {\sl $H$-balanceator} of $A$. More precisely, we have
\begin{eqnarray*}
\mathfrak{B}(a, b)&=&f^{a, b_2}\o g^{a, b_2}\o_H(t_{a, b_2}\o b_1)-g^{b_1, a}\o f^{b_1, a}\o_H (b_2\o t_{b_1, a})\\
&&+g^{b, a_2}\o f^{b, a_2}\o_H (a_1\o t_{b, a_2})-f^{a_1, b}\o g^{a_1, b}\o_H(t_{a_1, b}\o a_2),
\end{eqnarray*}
if $a*b=f^{a, b}\o g^{a, b}\o_H t_{a, b}$ for all $a, b\in A$.

The $H$-balanceator $\mathfrak{B}$ is said to be symmetric if $\mathfrak{B}(a, b)=(\sigma\o_H id)\mathfrak{B}(b, a)$ for all $a, b\in A$.

{\bf 5.1 From infinitesimal $H$-pseudobialgebras to Lie $H$-pseudobialgebras}.
\\

Recall that, if $(A,*)$ is an associative $H$-pseudoalgebra, then $(A, [*]_{lie})$ is a Lie $H$-pseudoalgebra obtained from the associative one, where $[a*b]_{lie}=a*b-(\sigma\o_H id)(b*a)$ for all $a, b\in A$. And if $(A, \D)$ is a coassociative $H$-coalgebra, then $(A,\delta_{lie}=(id-\tau)\D)$ is a Lie $H$-coalgebra by using the same argument as in the usual Lie coalgebra case. But, in general, we can not get that $(A,[*]_{lie},\delta_{lie})$ is a Lie $H$-pseudobialgebra due to $(A,*,\D)$ is an infinitesimal $H$-pseudobialgebra.
However, we have the following result.

{\bf Theorem 5.6} Let $(A, *, \D)$ be an infinitesimal $H$-pseudobialgebra. Then we have
\begin{eqnarray}
\nonumber&&\d_{lie}([a*b]_{lie})=[a*\d_{lie}(b)]_{lie}-(\sigma\o_H id)[b*\d_{lie}(a)]_{lie}+\mathfrak{B}(a, b)\\
&&-(\sigma\o_H id)\mathfrak{B}(b, a)
\end{eqnarray}
for all $a, b\in A$.

{\bf Proof.} We use the following shorthand: if $f(x, y)$ is an expression involving $x$ and $y$, then
\begin{eqnarray*}
&&(f(x, y))^{\wedge}=f(x, y)-(\sigma\o id)f(y, x).
\end{eqnarray*}
In particular, the right-hand side of equation $(5.3)$ becomes $([a*_{lie}\d(b)]_{lie}+\mathfrak{B}(a, b))^{\wedge}$.

Suppose $a*b=f^{a, b}\o g^{a, b}\o_H t_{a, b}$ for all $a, b\in A$. On the one hand, we have
\begin{eqnarray}
\nonumber&&\d_{lie}([a*b]_{lie})=(id-\tau)\D(a*b-(\sigma\o id)(b* a))\\
\nonumber&=&(id-\tau)\D(a*b)-(id-\tau)\D((\sigma\o id)(b* a))\\
\nonumber&=&f^{a, b_1}\o g^{a, b_1}\o_H(t_{a,b_1}\o b_2)+f^{a_2, b}\o g^{a_2, b}\o_H(a_1\o t_{a_2, b})\\
\nonumber&&-f^{a, b_1}\o g^{a, b_1}\o_H(b_2\o t_{a, b_1})-f^{a_2, b}\o g^{a_2, b}\o_H(t_{a_2, b}\o a_1)\\
\nonumber&&-g^{b, a_1}\o f^{b, a_1}\o_H(t_{b, a_1}\o a_2)-g^{b_2, a}\o f^{b_2, a}\o_H(b_1\o t_{b_2, a})\\
\nonumber&&+g^{b, a_1}\o f^{b, a_1}\o_H(a_2\o t_{b, a_1})+g^{b_2, a}\o f^{b_2, a}\o_H(t_{b_2, a}\o b_1)\\
\nonumber&=&(f^{a, b_1}\o g^{a, b_1}\o_H(t_{a,b_1}\o b_2)-g^{b_2, a}\o f^{b_2, a}\o_H(b_1\o t_{b_2£¬ a})\\
&&+g^{b_2, a}\o f^{b_2, a}\o_H(t_{b_2, a}\o b_1)-f^{a, b_1}\o g^{a, b_1}\o_H(b_2\o t_{a, b_1}))^{\wedge}.
\end{eqnarray}
On the other hand, observe that $A\o A$ is an $A$-module (here we view $A$ as a Lie $H$-pseudoalgebra with the pseudobracket $[*]_{lie}$) and the action is defined as in Lemma 4.1 in \cite{BL}. So we have
\begin{eqnarray}
\nonumber&&[a*\d_{lie}(b)]_{lie}=[a*(b_1\o b_2-b_2\o b_1)]\\
\nonumber&=&[a*b_1]\o b_2+b_1\o[a*b_2]-[a*b_2]\o b_1-b_2\o[a*b_1]\\
\nonumber&=&a*b_1\o b_2-(\sigma\o id)((b_1* a)\o b_2)+b_1\o(a* b_2)-b_1\o(\sigma\o id)(b_2* a)\\
\nonumber&&-a*b_2\o b_1+(\sigma\o id)((b_2* a)\o b_1)-b_2\o a*b_1+b_2\o(\sigma\o id)(b_1* a)\\
\nonumber&=&f^{a, b_1}\o g^{a, b_1}\o_H(t_{a, b_1}\o b_2)-g^{b_1, a}\o f^{b_1, a}\o_H(t_{b_1, a}\o b_2)+f^{a, b_2}\o g^{a, b_2}\o_H(b_1\o t_{a, b_2})\\
\nonumber&&-g^{b_2, a}\o f^{b_2, a}\o_H(b_1\o t_{b_2, a})-f^{a, b_2}\o g^{a, b_2}\o_H(t_{a, b_2}\o b_1)+g^{b_2, a}\o f^{b_2, a}\o_H(t_{b_2, a}\o b_1)\\
\nonumber&&-f^{a, b_1}\o g^{a, b_1}\o_H(b_2\o t_{a, b_1})+g^{b_1, a}\o f^{b_1, a}\o_H(b_2\o t_{b_1, a})\\
\nonumber&=&f^{a, b_1}\o g^{a, b_1}\o_H(t_{a, b_1}\o b_2)-g^{b_2, a}\o f^{b_2, a}\o_H(b_1\o t_{b_2, a})+g^{b_2, a}\o f^{b_2, a}\o_H(t_{b_2, a}\o b_1)\\
&&-f^{a, b_1}\o g^{a, b_1}\o_H(b_2\o t_{a, b_1})+(id\o_H \tau)([a, \D^{op}(b)]_{\star})-[a, \D^{op}(b)]_{\star}
\end{eqnarray}
Using equations $(5.4)$ and $(5.5)$, we obtain
\begin{eqnarray*}
&&([a*\d_{lie}(b)]_{lie}+\mathfrak{B}(a, b))^{\wedge}\\
&=&(f^{a, b_1}\o g^{a, b_1}\o_H(t_{a, b_1}\o b_2)-g^{b_2, a}\o f^{b_2, a}\o_H(b_1\o t_{b_2, a})+g^{b_2, a}\o f^{b_2, a}\o_H (t_{b_2, a}\o b_1)\\
&&-f^{a, b_1}\o g^{a, b_1}\o_H(b_2\o t_{a, b_1}))^{\wedge}+((id\o_H \tau)([a, \D^{op}(b)]_{\star})-[a, \D^{op}(b)]_{\star}\\
&&+[a, \D^{op}(b)]_{\star}+(\sigma\o_H \tau)([b, \D^{op}(a)]_{\star}))^{\wedge}\\
&=&(f^{a, b_1}\o g^{a, b_1}\o_H(t_{a, b_1}\o b_2)-g^{b_2, a}\o f^{b_2, a}\o_H(b_1\o t_{b_2, a})+g^{b_2, a}\o f^{b_2, a}\o_H (t_{b_2, a}\o b_1)\\
&&-f^{a, b_1}\o g^{a, b_1}\o_H(b_2\o t_{a, b_1}))^{\wedge}+(id\o_H \tau)([a, \D^{op}(b)]_{\star})-(\sigma\o_H \tau)([b, \D^{op}(a)]_{\star})\\
&&+(\sigma\o_H \tau)([b, \D^{op}(a)]_{\star})-(\sigma\o_H id)(\sigma\o_H \tau)([a, \D^{op}(b)]_{\star})\\
&=&\d_{lie}([a*b]_{lie}),
\end{eqnarray*}
which completes the proof.
$\hfill \square$
\\

The following result can be obtained directly by Theorem 5.6.

{\bf Corollary 5.7} Let $(A, *, \D)$ be an infinitesimal $H$-pseudobialgebra. Then $(A, [*]_{lie}, \d_{lie})$ is a Lie $H$-pseudobialgebra if and only if the $H$-balanceator $\mathfrak{B}$ is symmetric.
\\

Next, we discuss the construction of infinitesimal $H$-pseudobialgebras whose $H$-balancea-\\
tors are symmetric, then we can get related Lie $H$-pseudobialgebras
by Corollary 5.7.
\\

{\bf Proposition 5.8} Let $(A, *, \D)$ be an infinitesimal $H$-pseudobialgebra that is both commutative and cocommutative. Then its $H$-balanceator $\mathfrak{B}=0$.

{\bf Proof.} For all $a, b\in A$, we denote $a*b=(f^{a, b}\o g^{a, b})\o_H t_{a, b}$. Using the commutativity and cocommutativity of $A$, we compute:
\begin{eqnarray*}
&&f^{a_1, b}\o g^{a_1, b}\o_H(t_{a_1, b}\o a_2)+g^{b_1, a}\o f^{b_1, a}\o_H (b_2\o t_{b_1, a})\\
&=&g^{b, a_1}\o f^{b, a_1}\o_H(t_{b, a_1}\o a_2)+g^{b_2, a}\o f^{b_2, a}\o_H (b_1\o t_{b_2, a})\\
&=&(\sigma\o id)\D(b*a)=\D(a*b)\\
&=&f^{a_2, b}\o g^{a_2, b}\o_H(a_1\o t_{a_2, b})+f^{a, b_1}\o g^{a, b_1}\o_H (t_{a, b_1}\o b_2)\\
&=&g^{b, a_2}\o f^{b, a_2}\o_H(a_1\o t_{b, a_2})+f^{a, b_2}\o g^{a, b_2}\o_H (t_{a, b_2}\o b_1).
\end{eqnarray*}
So we have
\begin{eqnarray*}
\mathfrak{B}(a, b)&=&f^{a, b_2}\o g^{a, b_2}\o_H (t_{a, b_2}\o b_1)-g^{b_1, a}\o f^{b_1, a}\o_H (b_2\o t_{b_1, a})\\
&&+g^{b, a_2}\o f^{b, a_2}\o_H(a_1\o t_{b, a_2})-f^{a_1, b}\o g^{a_1, b}\o_H(t_{a_1, b}\o a_2)\\
&=&0,
\end{eqnarray*}
as required.
$\hfill \square$

{\bf 5.2 From coboundary infinitesimal $H$-pseudobialgebras to coboundary Lie $H$-pseudobialgebras}.
\\

In this subsection, we give a sufficient condition under which a coboundary infinitesimal $H$-pseudobialgebra gives rise to a coboundary Lie $H$-pseudobialgebra.

 Recall that a 2-tensor $r\in A\o A$ is said to be symmetric (resp, anti-symmetric) if $r=r^{op}$ (resp, $r=-r^{op}$), where $r^{op}=\tau(r)$. We start with a useful proposition in the following.
\\

{\bf Proposition 5.9} Let $(A, *, \D_r, r)$ be a coboundary infinitesimal $H$-pseudobialgebra with $r$ anti-symmetric. Then the $H$-balanceator $\mathfrak{B}$ of $A$ is $0$.

{\bf Proof.} For all $b\in A$, we have
\begin{eqnarray*}
&&\D_r(b)=h^{b, u_i}\c(c_{b, u_i}\o v_i)-l^{v_i, b}\c(u_i\o e_{v_i, b}),
\end{eqnarray*}
where $b* u_i=h^{b, u_i}\o 1\o_H c_{b, u_i}$ and $v_i* b=1\o l^{v_i, b}\o_H e_{v_i, b}$. It follows that
\begin{eqnarray}
\nonumber&&[a, \D_{r}^{op}(b)]_{\star}=a*\D_{r}^{op}(b)-(\sigma\o id)(\D_{r}^{op}(b)*a)\\
\nonumber&=&a*(h^{b, u_i}\c(v_i\o c_{b, u_i}))-a*(l^{v_i, b}\c(e_{v_i, b}\o u_i))-(\sigma\o id)(h^{b, u_i}\c(v_i\o c_{b, u_i})*a)\\
\nonumber&&+(\sigma\o id)(l^{v_i, b}\c(e_{v_i, b}\o u_i)* a)\\
\nonumber&=&(1\o h^{b, u_i}\o_H 1)(a*(v_i\o c_{b, u_i}))-(1\o l^{v_i, b}\o_H 1)(a*(e_{v_i, b}\o u_i))\\
\nonumber&&-(\sigma\o id)((h^{b, u_i}\o 1\o_H 1)((v_i\o c_{b, u_i})*a))+(\sigma\o id)((l^{v_i, b}\o 1\o_H 1)((e_{v_i, b}\o u_i)* a))\\
\nonumber&=&h^{a, v_i}\o h^{b, u_i}\o_H(c_{a, v_i}\o c_{b, u_i})-h^{a, e_{v_i, b}}\o l^{v_i, b}\o_H(c_{a, e_{v_i, b}}\o u_i)\\
&&-l^{c_{b, u_i}, a}\o h^{b, u_i}\o_H(v_i\o e_{c_{b, u_i}, a})+l^{u_i, a}\o l^{v_i, b}\o_H(e_{v_i, b}\o e_{u_i, a}).
\end{eqnarray}
Interchanging the roles of $a$ and $b$ in the above equation, we obtain
\begin{eqnarray}
\nonumber&&[b, \D_{r}^{op}(a)]_{\star}=h^{b, v_i}\o h^{a, u_i}\o_H(c_{b, v_i}\o c_{a, u_i})-h^{b, e_{v_i, a}}\o g^{v_i, a}\o_H (c_{b, e_{v_i, a}}\o u_i)\\
&&-g^{c_{a, u_i}, b}\o h^{a, u_i}\o_H(v_i\o e_{c_{a, u_i}, b})+g^{u_i, b}\o g^{v_i, a}\o_H(e_{v_i, a}\o e_{u_i, b}).
\end{eqnarray}
Using equations $(5.6)$ and $(5.7)$, we compute:
\begin{eqnarray}
\nonumber\mathfrak{B}(a, b)&=&[a, \D_r^{op}(b)]_{\star}+(\sigma\o \tau)([b, \D_r^{op}(a)]_{\star})\\
\nonumber&=&h^{a, v_i}\o h^{b, u_i}\o_H(c_{a, v_i}\o c_{b, u_i})-h^{a, e_{v_i, b}}\o l^{v_i, b}\o_H (c_{a, e_{v_i, b}}\o u_i)\\
\nonumber&&-l^{c_{b, u_i}, a}\o h^{b, u_i}\o_H(v_i\o e_{c_{b, u_i}, a})+l^{u_i, a}\o l^{v_i, b}\o_H(e_{v_i, b}\o e_{u_i, a})\\
\nonumber&&+h^{a, u_i}\o h^{b, v_i}\o_H(c_{a, u_i}\o c_{b, v_i})-l^{v_i, a}\o h^{b, e_{v_i, a}}\o_H (u_i\o c_{b, e_{v_i, a}})\\
\nonumber&&-h^{a, u_i}\o l^{c_{a, u_i}, b}\o_H(e_{c_{a, u_i}, b}\o v_i)+l^{v_i, a}\o l^{u_i, b}\o_H(e_{u_i, b}\o e_{v_i, a})\\
\nonumber&=&-h^{a, u_i}\o h^{b, v_i}\o_H(c_{a, u_i}\o c_{b, v_i})+h^{a, e_{u_i, b}}\o l^{u_i, b}\o_H (c_{a, e_{u_i, b}}\o v_i)\\
\nonumber&&+l^{c_{b, v_i}, a}\o h^{b, v_i}\o_H(u_i\o e_{c_{b, v_i}, a})-l^{v_i, a}\o l^{u_i, b}\o_H(e_{u_i, b}\o e_{v_i, a})\\
\nonumber&&+h^{a, u_i}\o h^{b, v_i}\o_H(c_{a, u_i}\o c_{b, v_i})-l^{v_i, a}\o h^{b, e_{v_i, a}}\o_H (u_i\o c_{b, e_{v_i, a}})\\
\nonumber&&-h^{a, u_i}\o l^{c_{a, u_i}, b}\o_H(e_{c_{a, u_i}, b}\o v_i)+l^{v_i, a}\o l^{u_i, b}\o_H(e_{u_i, b}\o e_{v_i, a})\\
\nonumber&=&h^{a, e_{u_i, b}}\o l^{u_i, b}\o_H (c_{a, e_{u_i, b}}\o v_i)-h^{a, u_i}\o l^{c_{a, u_i}, b}\o_H(e_{c_{a, u_i}, b}\o v_i)\\
&&+l^{c_{b, v_i}, a}\o h^{b, v_i}\o_H(u_i\o e_{c_{b, v_i}, a})-l^{v_i, a}\o h^{b, e_{v_i, a}}\o_H (u_i\o c_{b, e_{v_i, a}}).
\end{eqnarray}
Using the associativity of $A$, we have $a*(u_i*b)=(a*u_i)*b$ and $(b* v_i)*a=b*(v_i*a)$, which are equivalent to
\begin{eqnarray}
&&h^{a, m_{u_i, b}}\o 1\o g^{u_i, b}\o_H c_{a, m_{u_i, b}}=h^{a, u_i}\o 1\o g^{c_{a, u_i}, b}\o_H m_{c_{a, u_i}, b}
\end{eqnarray}
and
\begin{eqnarray}
&&h^{b, v_i}\o 1\o g^{c_{b, v_i}, a}\o_H m_{c_{b, v_i}, a}=h^{b, m_{v_i, a}}\o 1\o g^{v_i, a}\o_H c_{b, m_{v_i, a}}.
\end{eqnarray}
Now, comparing $(5.9)$ and $(5.10)$, all the terms on the right-hand side of $(5.8)$ are canceled. Thus $\mathfrak{B}=0$ as required.
$\hfill \square$
\\

{\bf Theorem 5.10} Let $(A, *, \D_r, r)$ be a coboundary infinitesimal $H$-pseudobialgebra with $r$ anti-symmetric. Then $A^{-}=(A, [*], \d_r, r)$ is a coboundary Lie $H$-pseudobialgebra, where $[a*b]=a*b-(\sigma\o id)(b* a), \d_r(a)=(id-\tau)\D_r(a)$ for all $a, b\in A$.

{\bf Proof.} By Proposition 5.9, the $H$-balanceator $\mathfrak{B}$ of $A$ is $0$, which is trivially symmetric. Thus Corollary 5.7 implies that $A^{-}$ is a Lie $H$-pseudobialgebra. It remains to show that $\d_r$ satisfies condition $(5.2)$. Set $r=\sum_{i} u_i\o v_i$ and $a*b=(f^{a, b}\o g^{a, b})\o_H t_{a, b}$ for all $a, b\in A$. On the one hand, using condition $r=-r^{op}$, we have
\begin{eqnarray}
\nonumber&&\d_r(a)=(id-\tau)\D_r(a)\\
\nonumber&=&(id-\tau)(f^{a, u_i}S(g^{a, u_i}_1)\c(g_2^{a, u_i}t_{a, u_i}\o v_i)-g^{v_i, a}S(f_1^{v_i, a})\c(u_i\o f_2^{v_i, a}t_{v_i, a}))\\
\nonumber&=&f^{a, u_i}S(g_{1}^{a, u_i})\c(g_{2}^{a, u_i}t_{a, u_i}\o v_i)-g^{v_i, a}S(f_{1}^{v_i, a})\c(u_i\o f_2^{v_i, a}t_{v_i, a})\\
\nonumber&&-f^{a, u_i}S(g_{1}^{a, u_i})\c(v_i\o g_{2}^{a, u_i}t_{a, u_i})+g^{v_i, a}S(f_{1}^{v_i, a})\c(f_2^{v_i, a}t_{v_i, a}\o u_i)\\
\nonumber&=&f^{a, u_i}S(g_{1}^{a, u_i})\c(g_{2}^{a, u_i}t_{a, u_i}\o v_i)-g^{u_i, a}S(f_{1}^{u_i, a})\c(f_2^{u_i, a}t_{u_i, a}\o v_i)\\
&&+f^{a, v_i}S(g_{1}^{a, v_i})\c(u_i\o g_{2}^{a, v_i}t_{a, v_i})-g^{v_i, a}S(f_{1}^{v_i, a})\c(u_i\o f_2^{v_i, a}t_{v_i, a}).
\end{eqnarray}
On the other hand, since
\begin{eqnarray*}
&&[a* u_i]=a*u_i-(\sigma\o id)(u_i* a)\\
&=&h^{a, u_i}\o g^{a, u_i}\o_H e_{a, u_i}-g^{u_i, a}\o h^{u_i, a}\o_H e_{u_i, a}\\
&=&h^{a, u_i}S(g_1^{a, u_i})\o 1\o_H g_2^{a, u_i}e_{a, u_i}-g^{u_i, a}S(h_1^{u_i, a})\o 1\o_H h_2^{u_i, a}e_{u_i, a},
\end{eqnarray*}
we have
$$
[a, u_i]=f^{a, u_i}S(g_1^{a, u_i})\o g_2^{a, u_i}t_{a, u_i}-g^{u_i, a}S(f_1^{u_i, a})\o f_2^{u_i, a}t_{u_i, a}.
$$
Similarly, we get
$$
[a, v_i]=f^{a, v_i}S(g_1^{a, v_i})\o g_2^{a, v_i}t_{a, v_i}-g^{v_i, a}S(f_1^{v_i, a})\o f_2^{v_i, a}t_{v_i, a}.
$$
Now we compute:
\begin{eqnarray}
\nonumber&&\mu([a, u_i]\o v_i+\sigma_{12}(u_i\o [a, v_i]))\\
\nonumber&=&(f^{a, u_i}S(g_1^{a, u_i}))\c(g_2^{a, u_i}t_{a, u_i}\o v_i)-(g^{u_i, a}S(f_1^{u_i, a}))\c(f_2^{u_i, a}t_{u_i, a}\o v_i)\\
&&+(f^{a, v_i}S(g_1^{a, v_i}))\c(u_i\o g_2^{a, v_i}t_{a, v_i})-(g^{v_i, a}S(f_1^{v_i, a}))\c(u_i\o f_2^{v_i, a}t_{v_i, a}).
\end{eqnarray}
Combining $(5.11)$ with $(5.12)$, we get $\d_r(a)=\mu([a, u_i]\o v_i+\sigma_{12}(u_i\o [a, v_i]))$, as desired.
$\hfill \square$

\section{From associative pseudo-Yang-Baxter equation to the classical type}\def\theequation{6.\arabic{equation}}
\setcounter{equation} {0}

Let $(L, [*])$ be a Lie $H$-pseudoalgebra and $r=\sum_{i} u_i\o v_i\in L\o L$. Recall that the classical Yang-Baxter equation ($pseudo$-CYBE) on $L$ has the form
\begin{eqnarray*}
&&[[r, r]]=\mu_{-1}^{3}([u_j, u_i]\o v_j\o v_i)-\mu_{-2}^{4}(u_i\o [u_j, v_i]\o v_j)-\mu_{-3}^{2}(u_i\o u_j\o[v_j, v_i])=0,
\end{eqnarray*}
where $\mu_{-k}^{l}$ means that the element of $H$ in the $k$-th place in its argument acts via the antipode on the element of $L$ located in the $l$-th place.
\\

{\bf Theorem 6.1} Let $(A, *)$ be an associative $H$-pseudoalgebra and $r=\sum_{i} u_i\o v_i\in A\o A$ a solution of the $pseudo$-AYBE. Suppose that $r$ is either symmetric or anti-symmetric. Then $r$ is a solution of the $pseudo$-CYBE in the Lie $H$-pseudoalgebra $A_{lie}=(A, [*]_{lie})$ for all $x, y\in A$.

{\bf Proof.} For all $a, b\in A$, we write $a*b=(f^{a, b}\o g^{a, b})\o_H t_{a, b}$. Then $[a*b]_{lie}=a*b-(\sigma\o_H id)(b*a)=(f^{a, b}\o g^{a, b})\o_H t_{a, b}-(g^{b, a}\o f^{b, a})\o_H t_{b, a}$. Define
\begin{eqnarray*}
&&A(r)'=\mu_{1}^{2,4}(\{u_j, u_i\}\o v_j\o v_i)-\mu_{2}^{1,3}(u_j\o \{u_i, v_j\}\o v_i)+\mu_{-3}^{1}(u_i\o u_j\o \{v_i, v_j\}),
\end{eqnarray*}
that is,
\begin{eqnarray*}
A(r)'&=& f_{1}^{u_j, u_i}t_{u_j, u_i}\o v_j\o f_2^{u_j, u_i}S(g^{u_j, u_i})v_i-f_2^{u_i, v_j}S(g^{u_i, v_j})u_j\o f_1^{u_i, v_j}t_{u_i, v_j}\o v_i\\
&&+g_1^{v_i, v_j}S(f^{v_i, v_j})u_i\o u_j\o g_2^{v_i, v_j}t_{v_i, v_j}.
\end{eqnarray*}
We first check that $A(r)'=\sigma_{13}(A(r))$. Using the (anti-)symmetry of $r$, we have
\begin{eqnarray*}
\sigma_{13}(A(r))&=&\sigma_{13}(\mu_{-1}^{4}(\{u_i, u_j\}\o v_j\o v_i)-\mu_{2}^{3,4}(u_i\o \{v_i, u_j\}\o v_j)\\
&&+\mu_{3}^{1,4}(u_i\o u_j\o \{v_j, v_i\}))\\
&=&\sigma_{13}(g_2^{u_i, u_j}t_{u_i, u_j}\o v_j\o g_1^{u_i, u_j}S(f^{u_i, u_j})v_i-u_i\o f_1^{v_i, u_j}t_{v_i, u_j}\o f_2^{v_i, u_j}S(g^{v_i, u_j})v_j\\
&&+f_1^{v_j, v_i}S(g^{v_j, v_i})u_i\o u_j\o f_2^{v_j, v_i}t_{v_j, v_i})\\
&=&g_1^{u_i, u_j}S(f^{u_i, u_j})v_i\o v_j\o g_2^{u_i, u_j}t_{u_i, u_j}-f_2^{v_i, u_j}S(g^{v_i, u_j})v_j\o f_1^{v_i, u_j}t_{v_i, u_j}\o u_i\\
&&+f_2^{v_j, v_i}t_{v_j, v_i}\o u_j\o f_1^{v_j, v_i}S(g^{v_j, v_i})u_i\\
&=&g_1^{v_i, v_j}S(f^{v_i, v_j})u_i\o u_j\o g_2^{v_i, v_j}t_{v_i, v_j}-f_2^{u_i, v_j}S(g^{u_i, v_j})u_j\o f_1^{u_i, v_j}t_{u_i, v_j}\o v_i\\
&&+f_2^{u_j, u_i}t_{u_j, u_i}\o v_j\o f_1^{u_j, u_i}S(g^{u_j, u_i})v_i\\
&=&A(r)'.
\end{eqnarray*}
Now, we compute:
\begin{eqnarray*}
&&[[r, r]]=\mu_{-1}^{3}([u_j, u_i]\o v_j\o v_i)-\mu_{-2}^{4}(u_i\o [u_j, v_i]\o v_j)-\mu_{-3}^{2}(u_i\o u_j\o[v_j, v_i])\\
&=&\mu_{-1}^{3}(f^{u_j, u_i}S(g_1^{u_j, u_i})\o g_2^{u_j, u_i}t_{u_j, u_i}\o v_j\o v_i-g^{u_i, u_j}S(f_1^{u_i, u_j})\o f_2^{u_i, u_j}t_{u_i, u_j}\o v_j\o v_i)\\
&&+\mu_{-2}^{4}(u_i\o g^{v_i, u_j}S(f_1^{v_i, u_j})\o f_2^{v_i, u_j}t_{v_i, u_j}\o v_j-u_i\o f^{u_j, v_i}S(g_1^{u_j, v_i})\o g_2^{u_j, v_i}t_{u_j, v_i}\o v_j)\\
&&+\mu_{-3}^{2}(u_i\o u_j\o g^{v_i, v_j}S(f_1^{v_i, v_j})\o f_2^{v_i, v_j}t_{v_i, v_j}-u_i\o u_j\o f^{v_j, v_i}S(g_1^{v_j, v_i})\o g_2^{v_j, v_i}t_{v_j, v_i})\\
&=&g_2^{u_j, u_i}t_{u_j, u_i}\o g_1^{u_j, u_i}S(f^{u_j, u_i})v_j\o v_i-f_2^{u_i, u_j}t_{u_i, u_j}\o f_1^{u_i, u_j}S(g^{u_i, u_j})v_j\o v_i\\
&&+u_i\o f_2^{v_i, u_j}t_{v_i, u_j}\o f_1^{v_i, u_j}S(g^{v_i, u_j})v_j-u_i\o g_2^{u_j, v_i}t_{u_j, v_i}\o g_1^{u_j, v_i}S(f^{u_j, v_i})v_j\\
&&+u_i\o f_1^{v_i, v_j}S(g^{v_i, v_j})u_j\o f_2^{v_i, v_j}t_{v_i, v_j}-u_i\o g_1^{v_j, v_i}S(f^{v_j, v_i})u_j\o g_2^{v_j, v_i}t_{v_j, v_i}\\
&=&f_2^{u_j, u_i}t_{u_j, u_i}\o v_j\o f_1^{u_j, u_i}S(g^{u_j, u_i})v_i-g_2^{u_i, u_j}t_{u_i, u_j}\o v_j\o g_1^{u_i, u_j}S(f^{u_i, u_j})v_i\\
&&+u_i\o f_1^{v_i, u_j}t_{v_i, u_j}\o f_2^{v_i, u_j}S(g^{v_i, u_j})v_j-f_2^{u_i, v_j}S(g^{u_i, v_j})u_j\o f_1^{u_i, v_j}t_{u_i, v_j}\o v_i\\
&&+ g_1^{v_i, v_j}S(f^{v_i, v_j})u_i\o u_j\o g_2^{v_i, v_j}t_{v_i, v_j}-f_1^{v_j, v_i}S(g^{v_j, v_i})u_i\o u_j\o f_2^{v_j, v_i}t_{v_j, v_i}\\
&=&A(r)'-A(r)\\
&=&(\sigma_{13}-id)(A(r))\\
&=&0.
\end{eqnarray*}
So the conclusion holds.
$\hfill \square$
\\

{\bf Example 6.2} Consider the free associative $H$-pseudoalgebra $A=H\{e_1, e_2\}$ in Example 4.3 with pseudoproduct
$$
e_1*e_1=e_1*e_2=0, \quad e_2*e_1=1\o 1\o_H e_1, \quad e_2*e_2=1\o 1\o_H e_2.
$$
Then $r=e_2\o e_1-e_1\o e_2$ is an anti-symmetric solution of $pseudo$-AYBE. By Theorem 6.1, $r$ induces a solution of $pseudo$-CYBE in $A_{lie}$, where the pseudobracket in $A_{lie}$ is determined by $[e_1*e_1]_{lie}=[e_2*e_1]_{lie}=[e_2*e_2]_{lie}=0, [e_1*e_2]_{lie}=-1\o 1\o_H e_1$.

 {\bf Acknowledgments.}
 This work was partially supported by the Doctoral Foundation of Henan Institute of Technology (grant no. KQ2003).

 \end{document}